\input amstex
\magnification=\magstep1
\documentstyle{amsppt}
\advance\hoffset by 0.1 truein
\TagsOnLeft
\hoffset=0.23in
\hsize=4.87in
\voffset=0in
\vsize=7.5in
\parindent=12pt
\font\smallit=cmmi6
\font\smallrm=cmr8
\font\smallslant=cmsl8
\document

\centerline{\bf COMPACTIFYING THE RELATIVE JACOBIAN}

\vskip0.1cm

\centerline{\bf OVER FAMILIES OF REDUCED CURVES}

\vskip0.5cm

\centerline{\smc eduardo esteves\footnote 
{Research supported by an MIT Japan Program 
Starr fellowship and by CNPq, Proc. 300004/95-8.}}

\vskip0.5cm

\author 
eduardo esteves
\endauthor

\title 
compactifying the relative jacobian
\endtitle

\vskip1.0cm

\heading 
Contents
\endheading

\parindent=0pt

0. Introduction

1. Semistable sheaves

2. The valuative criteria

3. The fine moduli spaces

4. Certain special cases

5. Theta functions and the separated realizations

6. Comparison with Seshadri's compactification

References

\vskip0.8cm

\heading 
0. Introduction
\endheading

\parindent=0pt

\subheading{0.1. Motivation and Goal} 
The problem of finding a natural relative 
compactification of the relative Jacobian over a family of curves 
has drawn a lot of attention since Igusa's pioneering work \cite{\bf 18} 
in the fifties. (In Subsect. 0.2 the reader will find a slightly 
more detailed account of 
this history, with references.) 
In the case where the curves in the family are geometrically 
integral, a very satisfactory solution has been found by 
Altman and Kleiman \cite{\bf 4}: their relative compactification 
is a fine moduli space; that is, it admits a 
universal object, after an \'etale base change. 
However, reducible curves (especially nodal ones) show up quite often in 
applications. For instance, Deligne-Mumford stable curves 
are used in the compactification, $\overline M_g$, 
of the moduli space of non-singular curves of genus $g$. 
Recently, Caporaso \cite{\bf 10} and Pandharipande 
\cite{\bf 25} produced a 
relative compactification of the relative Jacobian over $\overline M_g$. Their 
construction is strongly based on Gieseker's construction \cite{\bf 15} 
of $\overline M_g$, and does not seem to be adaptable to different 
situations. The main 
disadvantages of the constructions found so far for reducible curves 
are the lack of a 
universal object and the restricted range of applicability.

\parindent=12pt

Apart from the study of $\overline M_g$, the relatively compactified Jacobian 
has been most recently employed by Beauville \cite{\bf 8} in counting 
the number of rational curves on K3 surfaces. In his article, 
Beauville made the simplifying assumption \cite{\bf 8\rm , 1.2} 
that all curves belonging to a certain linear system on the K3 surface 
are integral. The assumption was used in order to guarantee the 
existence of a fine relative compactification of the relative Jacobian 
over the family of curves belonging to the linear system. 
In general, however, certain curves of the linear 
system may be reducible, and may have singularities more complicated than 
ordinary nodes. 

These recent developments suggest the urge for a more general approach 
to the problem of compactifying the relative Jacobian. 
One approach was taken recently by 
Simpson \cite{\bf 29\rm , Sect. 1}, 
who constructed moduli spaces of coherent sheaves on any family of 
varieties by means of Geometric Invariant Theory (G.I.T.). The main 
disadvantage of his moduli spaces is, again, the lack of a universal 
object.

Roughly speaking, the present article aims at constructing a natural relative 
compactification of the relative Jacobian over 
a projective, flat family of geometrically reduced and connected curves. 
In contrast to earlier relative compactifications, ours admits a universal 
object, after an \'etale base change. The points of our 
compactification correspond to simple, torsion-free, rank 1 sheaves 
that are semistable with respect to a given polarization. To compare 
our relative compactification with that obtained by Seshadri 
in \cite{\bf 27\rm , Part 7}, we use theta functions. 
In contrast to all the past approaches, we do not use 
G.I.T.. (The method of theta functions has been already 
used by Faltings to construct 
the moduli of semistable vector bundles on a 
smooth complete curve (see \cite{\bf 14} or 
\cite{\bf 28}).) It must be said that our relative compactification is 
an algebraic space, rather than a scheme. (It could not be 
otherwise, as a famous example by Mumford \cite{\bf 9\rm , p. 210} 
shows that the relative Jacobian is not a scheme, in general.) But we show 
that it becomes a scheme, after an \'etale base change. The reader will 
find a more detailed account of our results in Subsect. 0.3.

\vskip0.4cm

\parindent=0pt

\subheading{0.2. History} Igusa \cite{\bf 18} was probably the 
first to consider the problem of compactifying the (generalized) 
Jacobian variety of a reduced curve $X$. His method was to 
construct a compactification as the limit of the Jacobians of 
smooth curves approaching $X$, and he applied his method to the case where 
$X$ was nodal and irreducible. He showed also that his 
compactification, in spite of his construction, did not depend on the 
family of approaching smooth curves. Later, 
Meyer and Mumford \cite{\bf 22} announced an intrinsic characterization of 
Igusa's compactification by means of torsion-free, rank 1 sheaves. 
Their idea was carried out by D'Souza in 
his thesis \cite{\bf 11}, where he used G.I.T. to construct a 
compactification of the Jacobian variety, when $X$ is irreducible, 
and then showed that there is a universal torsion-free, rank 1 sheaf 
over it. In the case where $X$ is reducible and nodal, 
Oda and Seshadri \cite{\bf 24} used G.I.T. to construct several 
compactifications of disjoint unions of copies of the Jacobian variety. 
Finally, Seshadri used G.I.T. to deal with a general reduced curve $X$ 
in \cite{\bf 27} (where he considered also the higher rank case).

\parindent=12pt

In the case of families of curves, the (relative) 
compactification problem is more 
difficult, as the relative Jacobian itself 
does not behave well. Despite that, some results have been obtained, as I 
will mention below.

Let $S$ be a locally Noetherian scheme. 
Let $f\: X @>>> S$ be a projective, flat morphism whose geometric 
fibres are connected, reduced curves. Let
$$
\text{\bf P}^*\: (\text{Sch}/S)^o @>>> (\text{Sets})
$$
be the \sl relative Jacobian functor\rm , 
defined on an $S$-scheme $T$ as the set of invertible sheaves on 
$X\times_S T$. Let 
$\text{\bf P}$ be the \'etale sheaf associated to $\text{\bf P}^*$. 
Artin \cite{\bf 7} 
showed that the functor $\text{\bf P}$ is represented by an algebraic 
space $P$, locally 
of finite type over $S$. If the geometric fibres of $f$ are integral, 
then Grothendieck 
showed that $P$ is a scheme \cite{\bf 17}. 
Without any assumption on the fibres, Grothendieck showed also that, 
if $S$ is integral, then there is an open dense subscheme $S'\subseteq S$ 
such that $P\times_S S'$ is a scheme 
\cite{\bf 26}. 
In general, $P$ is not a scheme, as a famous 
example by Mumford shows \cite{\bf 9\rm , p. 210}. 
But Mumford himself 
proved that, if the irreducible components of the fibres of $f$ are 
geometrically integral, then $P$ is a 
scheme \cite{\bf 9\rm , p. 210}. 

The algebraic space $P$ is formally smooth over $S$, 
but may not be either separated or of finite type over $S$. 
If the geometric fibres of $f$ are integral, then the subspace 
$P_d\subseteq P$, parametrizing invertible sheaves with Euler 
characteristic $d$, is a 
separated scheme of finite type over $S$ \cite{\bf 17}. 
In general, however, $P_d$ may not even be of finite type over $S$.

The problem of compactifying $P_d$ over $S$ 
has attracted a lot of attention from 
the start. D'Souza himself \cite{\bf 11} worked out a compactification of 
$P_d$ over $S$, when $S$ is (the 
spectrum of) a Henselian local ring with separably closed 
residue field, and the geometric fibres of $f$ are integral. 
More generally though, Altman and Kleiman realized that 
Meyer's and Mumford's idea 
could be combined with Mumford's \cite{\bf 23} or Grothendieck's 
\cite{\bf 17} method of constructing the Picard scheme to yield a 
natural relative compactification of 
$P_d$ over $S$ by means of torsion-free, rank 1 sheaves, when the 
geometric fibres of 
$f$ are integral (see \cite{\bf 3} and \cite{\bf 4}, 
respectively, where the higher dimesion case is also considered). 
When the geometric fibres of $f$ 
are nodal (possibly reducible) curves, then 
Ishida \cite{\bf 19} adapted Oda's and 
Seshadri's method to describe several compactifications of subspaces 
of $P_d$, but his hypotheses are numerous. 

After a long time without any 
work on the relative compactification problem, Caporaso \cite{\bf 10} 
showed how to compactify the relative Jacobian over the moduli of stable 
curves by putting invertible sheaves on curves derived from stable curves on 
the boundary. One year 
later, Pandharipande \cite{\bf 25} 
produced the same compactification, with the boundary points now representing 
torsion-free, rank 1 sheaves, following 
Seshadri's method in \cite{\bf 27}. As we have already mentioned, 
both Caporaso's and Pandharipande's construction rely heavily on 
Gieseker's construction of 
the moduli of stable curves \cite{\bf 15}, and thus cannot be extended to 
an arbitrary family of curves.

Recently, Simpson \cite{\bf 29\rm , Sect. 1} constructed moduli spaces 
of coherent sheaves over any family of projective varieties.

\vskip0.8cm

\parindent=0pt

\subheading{0.3. The results} We present now a simplified account of the 
methods and contributions of the present article. 
Let $S$ be a locally Noetherian scheme. Let $f\: X @>>> S$ 
be a projective, flat morphism whose geometric fibres are connected, 
reduced curves. Let $P$ denote the algebraic space over $S$ 
parametrizing invertible sheaves on $X/S$ (see Subsect. 0.2). 
In order to compactify $P$ over $S$ in general, 
it would be natural, following Meyer and Mumford, to consider the functor
$$
\text{\bf F}^* \: (\text{Sch}/S)^o @>>> (\text{Sets}),
$$
defined on an $S$-scheme $T$ as the set of $T$-flat, coherent sheaves on 
$X\times_S T$ whose fibres over $T$ are 
torsion-free, rank 1 sheaves. Defining $\text{\bf F}$ 
as the \'etale sheafification of $\text{\bf F}^*$, it is clear that 
$\text{\bf F}$ contains $\text{\bf P}$ as an open subfunctor. 
It is also easy to show that 
$\text{\bf F}$ ``contains enough degenerations''. In other words, 
$\text{\bf F}$ meets the 
existence condition of the valuative criterion of properness, without 
necessarily meeting the uniqueness condition. 
However, the functor $\text{\bf F}$ 
is not representable by an algebraic space in 
general, the obstruction being the existence of 
torsion-free, rank 1 sheaves that are not simple (i.e. admit non-trivial 
endomorphisms). In fact, considering the subfunctor 
$\text{\bf J}\subseteq \text{\bf F}$, parametrizing sheaves with 
simple fibres, it 
follows from \cite{\bf 4\rm , Thm. 7.4, p. 99} 
that $\text{\bf J}$ is represented 
by an algebraic space $J$. It is clear that $J$ contains $P$ as 
an open subspace, since the geometric fibres of $f$ are reduced and 
connected. But does $J$ ``contain enough degenerations'' over $S$? 
The surprising answer we obtain in Thm. 25 of 
the present article is: yes! The upshot is that we do not need to 
consider all torsion-free, rank 1 
sheaves to compactify the relative Jacobian, but just those that are 
simple. Of course, $J$ is neither separated nor of finite type 
over $S$ in general, since neither is the open subspace $P$. Nevertheless, 
since $J$ is a fine moduli space, it is worthwhile to analyse $J$ and, 
perhaps, obtain from $J$ ``coarse moduli spaces'' that behave better than 
$J$. The present article is thus devoted to the study of $J$.

\parindent=12pt

Since $J$ is awkwardly ``big'', we need to decompose $J$ into simpler 
``pieces'', and for that we use polarizations like those defined by 
Seshadri in \cite{\bf 27\rm , Part 7, p. 153}. As a matter of fact, 
Seshadri used numerical 
polarizations but, since we want to deal with a family of curves, 
we prefer to use ``continuous'' polarizations. For us, a 
polarization on $X$ over $S$ will be a vector bundle $\Cal E$ on $X$ 
with rank $r>0$ and relative degree $-rd$ over $S$, for a certain integer $d$ 
(see Sect. 1). Of course, there are natural 
choices of relative 
polarizations, as the structure sheaf $\Cal O_X$ and those 
constructed from $\Cal O_X$ and the relative dualizing sheaf 
$\omega$, when the fibres of $f$ are 
Gorenstein. (The latter were used by Pandharipande \cite{\bf 25}.)

Associated to a polarization $\Cal E$ on $X$ over $S$ we have the usual 
classes of stable and semistable torsion-free, rank 1 sheaves. In 
Sect. 1 we define two new classes of sheaves, those of quasistable and 
$\sigma$-quasistable torsion-free, rank 1 sheaves (see 
Subsections 1.2 and 1.5), where 
$\sigma\: S @>>> X$ is a section of $f$ through the 
$S$-smooth locus of $X$. These new classes of 
sheaves are important for their cohomological and 
degeneration properties (see Subsect. 1.4 and Sect. 2). Let 
$J^s_{\Cal E}$ (resp. $J^{ss}_{\Cal E}$, resp. $J^{qs}_{\Cal E}$, resp. 
$J^{\sigma}_{\Cal E}$) be the 
subspace of $J$ parametrizing sheaves with 
stable (resp. semistable, resp. quasistable, resp. 
$\sigma$-quasistable fibres) with respect to the polarization $\Cal E$. 
It follows from the cohomological characterizations of Subsect. 1.4 that all 
the above subspaces are open in $J$. 
In general, neither $J_{\Cal E}^{qs}$ nor $J_{\Cal E}^{ss}$ is separated over 
$S$. These spaces are still ``too big''. Nevertheless, applying the 
method of Langton's \cite{\bf 21} (see Sect. 2), 
we prove the following theorem.

\vskip0.4cm

\proclaim{Theorem A} The algebraic space 
$J_{\Cal E}^{ss}$ is of finite type over $S$. In addition:
\roster
\item $J_{\Cal E}^{ss}$ and $J_{\Cal E}^{qs}$ are 
universally closed over $S$;
\item $J_{\Cal E}^s$ is separated over $S$;
\item $J_{\Cal E}^{\sigma}$ is proper over $S$.
\endroster
\endproclaim

\vskip0.4cm

Provided it is convenient 
to fix a section $\sigma$ of $f$ through the $S$-smooth locus of $X$ (the 
case of a family of pointed curves), we may restrict ourselves to 
$J_{\Cal E}^{\sigma}$. We can always make an \'etale base change to 
obtain enough sections of $f$. In fact, a suitable \'etale base change 
will also give us schemes, as it follows from the next theorem.

\vskip0.4cm

\proclaim{Theorem B} Assume that there are sections 
$\sigma_1,\dots,\sigma_n\: S @>>> X$ of $f$ such that:
\roster
\item $\sigma_i$ factors through the $S$-smooth locus of $X$ for 
$i=1,\dots,n$;
\item for every $s\in S$, every irreducible component of $X(s)$ is 
geometrically integral, and contains $\sigma_i(s)$ for some $i$.
\endroster 
Then, $J$ is a scheme.
\endproclaim

\vskip0.4cm

However, it might not be appropriate, sometimes, to 
change the base scheme, or even to fix a section $\sigma$ of 
$f$ through the $S$-smooth locus of $X$, assuming it exists. 
To overcome these problems, we can use 
theta functions (associated to vector bundles on $X$) 
to construct ``approximations'' of the 
algebraic spaces $J_{\Cal E}^{ss}$, $J_{\Cal E}^{qs}$ and  
$J_{\Cal E}^{\sigma}$ that are locally projective schemes over $S$. 
For simplicity, assume that 
$S$ is locally of finite type over an algebraically closed field $k$. 
Let $\Sigma\subseteq J_{\Cal E}^{ss}$ be an open subspace, and denote 
by $\pi\: \Sigma @>>> S$ the structure morphism. 
Then, there is a natural invertible sheaf 
$\Cal L_{\Cal E}$ on $\Sigma$, uniquely defined from the polarization 
$\Cal E$. Moreover, there is a natural quasicoherent graded subsheaf of 
$\Cal O_S$-subalgebras,
$$
V_{\Cal E}\subseteq \Gamma_{\Cal E}:=\bigoplus_{m\geq 0} 
\pi_*\Cal L_{\Cal E}^{\otimes m},
$$
whose homogeneous pieces are coherent. The sheaf $V_{\Cal E}$ is 
generated by the so-called theta functions (associated to vector 
bundles on $X$). Let $\overline{\Sigma}:=\text{Proj}(V_{\Cal E})$, and 
consider the natural rational map, $\psi\: \Sigma @>>> 
\overline{\Sigma}$.

\vskip0.4cm

\proclaim{Theorem C} Assume that $S$ is locally of finite type over 
an algebraically closed field $k$. 
Let $\Sigma\subseteq J_{\Cal E}^{ss}$ be an open 
subspace. Then, the following statements hold:
\roster
\item $\overline{\Sigma}$ is locally projective over $S$.
\item The rational map $\psi$ is defined on $\Sigma$, and is 
scheme-theoretically dominant.
\item For every closed point $s\in S$, each fibre of $\psi(s)$ 
is contained in a Jordan-H\"older equivalence class of $\Sigma(s)$.
\item If $\Sigma$ is universally closed over $S$, and 
$J^s_{\Cal E}\subseteq\Sigma$, then the restriction of $\psi$ to 
$J_{\Cal E}^s$ is an open embedding.
\endroster
\endproclaim

\vskip0.4cm

Finally, we compare our compactification with Seshadri's 
\cite{\bf 27\rm , Part 7}, 
in the case where $S$ is (the spectrum of) an algebraically 
closed field $k$ (see Sect. 6). Roughly speaking, we show that, if 
$X$ has at most ordinary double points for singularities, then 
any structure theta functions detect in Seshadri's compactification is also 
detected in ours, and vice-versa. 
(See Thm. 48, and the discussion before it, for a precise statement.) 

Some conventions: by a curve $X$ we shall mean a geometrically 
reduced scheme of pure dimension 1, and of finite type over a field $k$; 
by a vector bundle we mean a locally free sheaf of 
constant rank; all our schemes will be assumed locally Noetherian. 

\vskip0.8cm

\parindent=12pt

\heading 
1. Semistable sheaves
\endheading

\parindent=0pt

\subheading{1.1. Preliminaries} Let $X$ be a geometrically reduced curve 
over a field $k$. Let $X_1,\dots,X_n$ denote the irreducible components of 
$X$. In this section, we shall assume 
that $X_1,\dots,X_n$ are geometrically integral. 
A \sl subcurve \rm of $X$ is a reduced subscheme $Y\subseteq X$ of 
pure dimension 1. The empty set will also be considered a subcurve of $X$. 
If $Y,Z\subseteq X$ are 
subcurves, we let $Y\wedge Z$ denote the maximum subcurve of $X$ 
contained in $Y\cap Z$; and we let $Z-Y$ denote the minimum subcurve 
containing $Z\setminus Y$. If $Y\subseteq X$ is a subcurve, we let 
$Y^c:=X-Y$.

\parindent=12pt

Let $I$ be a coherent sheaf on $X$. We say that $I$ is \sl torsion-free \rm 
if $I$ has no embedded components. We say that $I$ has \sl rank 1 \rm if 
$I$ has generic rank 1 at every irreducible component of $X$. We 
say that $I$ is \sl simple \rm if $\text{End}_X(I)=k$. 
If $I$ is invertible, then $I$ is 
torsion-free, rank 1, and also simple if $X$ is connected.

Let $I$ be a torsion-free, rank 1 sheaf on $X$. If $Y\subseteq X$ is a
subcurve, then we will denote by $I_Y$ the maximum torsion-free quotient of 
$\left.I\right|_Y$. Of course, there is a canonical
surjective homomorphism, $I \twoheadrightarrow I_Y$. We say that 
$I$ is \sl decomposable \rm if there are proper subcurves 
$Y,Z\subsetneqq X$ such that the canonical homomorphism, 
$I @>>> I_Y\oplus I_Z$, is an isomorphism. In this case, we say that $I$ 
\sl decomposes \rm at $Y$ (or $Z$).

\vskip0.4cm

\proclaim{Proposition 1} Let $I$ be a torsion-free, rank 1
sheaf on $X$. Then $I$ is simple if and only if $I$ is not
decomposable. 
\endproclaim

\demo{Proof} It is clear that, if $I$ is decomposable,
then $I$ is not simple. Suppose that $I$ is not simple.
Then, there is an endomorphism $h\: I @>>> I$ that is
not a multiple of the identity. Let $Y\subseteq X$ be the subcurve such that 
$I_Y\cong\text{im}(h)$, and let $h'\: I_Y \hookrightarrow I$ 
denote the induced injective homomorphism. Since $h\neq 0$, then 
$Y$ is not empty. Since $I_W$ is simple for any irreducible 
component $W\subseteq X$, by \cite{\bf 4\rm , Lemma 5.4, p. 83}, 
then we may also assume that $Y\neq
X$. It is clear that $h'$ factors through 
$J:=\text{ker}(I \twoheadrightarrow I_Z)$, 
where $Z:=Y^c$. Moreover, since $h'$ and the 
composition, $J \hookrightarrow I \twoheadrightarrow I_Y$, are injective 
homomorphisms, then $h'$ is actually an isomorphism onto 
$J$. In particular, $I$ and $I_Y\oplus I_Z$ have the same
Euler characteristic. Thus, $I=I_Y\oplus I_Z$. The proof is complete.\qed
\enddemo

\vskip0.4cm

The above proposition does not hold in higher rank, even if we assume that 
$X$ is smooth. In fact, if $X$ is smooth and not rational, then any 
vector bundle $E$ fitting in the middle of a non-split short 
exact sequence of the form:
$$
0 @>>> \Cal O_X @>>> E @>>> \Cal O_X @>>> 0
$$
is neither simple nor decomposable. Though easy to state and prove, 
the above proposition is the key 
reason why we are able to obtain fine moduli spaces in the rank 1 
case.

\vskip0.4cm

\proclaim{Lemma 2} Let $Y,Z\subseteq X$ be 
non-empty subcurves covering $X$. 
Let $M$ be a torsion-free, rank 1 sheaf on $X$. Then, the 
following statements hold:
\roster 
\item If $Y\wedge Z\neq\emptyset$, and both $M_Y$ and $M_Z$ are simple, 
then $M$ is simple.
\item If there is an exact sequence of the 
form:
$$
0 @>>> I @>>> M @>>> J @>>> 0,
$$
where $I$ (resp. $J$) is a simple, torsion-free, rank 1 sheaf on 
$Y$ (resp. $Z$), then $M$ is simple if and only if the sequence 
is not split.
\endroster
\endproclaim

\demo{Proof} Suppose that there are 
subcurves $X_1,X_2\subseteq X$ such that $M=M_{X_1}\oplus 
M_{X_2}$. In addition, suppose that there is a surjective homomorphism, 
$\mu\: M \twoheadrightarrow J$, where 
$J$ is a simple, torsion-free, rank 1 sheaf on a subcurve $Z\subseteq X$. 
Of course, $\mu$ is the direct sum of two homomorphisms,
$$
\mu_1\: M_{X_1} @>>> J \text{\  \  and \  \  } 
\mu_2\: M_{X_2} @>>> J.
$$
Since $\mu$ is surjective, then $\text{im}(\mu_i)=M_{Z_i}$, where 
$Z_i:=Z\wedge X_i$ for $i=1,2$, and
$$
J=\text{im}(\mu)=\text{im}(\mu_1)\oplus \text{im}(\mu_2).
$$
Since $J$ is simple, then either $Z\subseteq X_1$ or $Z\subseteq X_2$.

We prove (1) now. Suppose that $M=M_{X_1}\oplus M_{X_2}$, for 
subcurves $X_1,X_2\subseteq X$. Apply the above reasoning twice, 
to both $J:=M_Y$ and $J:=M_Z$. Without loss of generality, we either have that 
$Y\subseteq X_1$ and $Z\subseteq X_2$, or $Y\cup Z\subseteq X_1$. Since 
$Y\wedge Z\neq\emptyset$, then $Y\cup Z\subseteq X_1$. Since $Y\cup Z=X$, then 
$X_1=X$. It follows from Prop. 1 that $M$ is simple.

We prove (2) now. The ($\Rightarrow$) part 
is trivial. We show ($\Leftarrow$) now. Suppose by contradiction that 
$M=M_{X_1}\oplus M_{X_2}$, for 
proper subcurves $X_1,X_2\subsetneqq X$. Apply the reasoning in the first 
paragraph of the proof 
to the surjection $\mu\: M \twoheadrightarrow J$. 
We may assume, without loss of 
generality, that $Z\subseteq X_1$. Then $\mu_2=0$, and hence 
$I=\text{ker}(\mu)=\text{ker}(\mu_1)\oplus M_{X_2}$. Since $I$ is simple, 
and $X_2\neq\emptyset$, then 
$\text{ker}(\mu_1)=0$. So $J=M_{X_1}$, and thus the sequence in (2) 
is split, a contradiction. The proof is complete.\qed
\enddemo

\vskip0.4cm

We observe that, if $I$ is a torsion-free, rank 1 sheaf on $X$, and 
$E$ is a vector bundle on $X$ of rank $r$ and degree $e$, then 
$\chi(I\otimes E)=r\chi(I)+e$.

\vskip0.4cm

\parindent=0pt

\subheading{1.2. Semistable sheaves} Fix an integer $d$. Fix a vector
bundle $E$ on $X$ of rank $r>0$ and degree $-rd$. 
We say that $E$ is a \sl polarization \rm on $X$. 
For every subcurve $Y\subseteq X$, let 
$e_Y:=-\deg_Y E$. 
Note that $E$ induces the polarization $\left. E\right|_Y$ on $Y$ as long as 
$r|e_Y$.

\parindent=12pt

Let $I$ be a torsion-free, rank 1 sheaf on $X$ with $\chi(I)=d$. 
We say that $I$ is \sl stable \rm (resp. \sl semistable\rm )
\sl with respect to $E$ \rm if, for every non-empty, proper 
subcurve $Y\subsetneqq X$,
$$
\chi(I_Y)>e_Y/r \text{ (resp. }\chi(I_Y)\geq e_Y/r
\text{)},
$$
or equivalently,
$$
\chi(I_Y\otimes E)>0 \text{ (resp. }\chi(I_Y\otimes E)\geq
0\text{)}.
$$
If $X$ is irreducible, then any torsion-free, rank 1 sheaf $I$ on $X$
with $\chi(I)=d$ is stable with respect to $E$. 

Strictly 
speaking, the notions of stability 
and semistability do
not depend on $E$, but rather on its \sl multi-slope\rm :
$$
\underline \mu_E:=(e_{X_1},\dots,e_{X_n})/r.
$$
However, in dealing with families of curves, we will find it more convenient 
to think of $E$ as the polarization, rather than $\underline \mu_E$. 
When the polarization $E$ is clear from the context, we will seldom make 
reference to $E$.

The above notions of stability and 
semistability coincide with those of \cite{\bf 13\rm , Sect. 3}. 
We refer to \cite{\bf 13\rm , Rmk. 12} for a comparison between our notions of 
polarization, semistability and stability 
with Seshadri's \cite{\bf 27\rm , Part 7}.

Let $I$ be a torsion-free, rank 1 sheaf on $X$ with $\chi(I)=d$. 
For every subcurve $Y\subseteq X$, let
$\beta_I(Y):=\chi(I_Y)-e_Y/r$. Of course, $I$ is stable (resp.
semistable) if and only if 
$\beta_I(Y)>0$ (resp. $\beta_I(Y)\geq 0$) 
for every non-empty, proper subcurve $Y\subsetneqq X$. 
If $I$ is semistable, and $Y\subseteq 
X$ is a subcurve, then 
$\beta_I(Y)=0$ if and only if
$I_Y$ is semistable with respect to $\left. E\right|_Y$.

\vskip0.4cm

\proclaim{Lemma 3} Let $I$ be a torsion-free, rank 1
sheaf on $X$ with $\chi(I)=d$. If $Y,Z\subseteq X$ are subcurves,
then
$$
\chi(I_{Y\cup Z})+\chi(I_{Y\wedge Z})\leq \chi(I_Y)+\chi(I_Z),
$$
or equivalently,
$$
\beta_I(Y\cup Z)+\beta_I(Y\wedge Z)\leq \beta_I(Y)+\beta_I(Z).
$$
\endproclaim

\demo{Proof} The inequalities follow immediately from
the following canonical commutative diagram:
$$
\CD
I_{Y\cup Z} @>>> I_Y\oplus I_{Z-Y}\\
@VVV @VVV\\
I_Z @>>> I_{Y\wedge Z}\oplus I_{Z-Y},
\endCD
$$
when we observe that the vertical homomorphisms are surjective, and the
horizontal ones are injective with finite length cokernel.\qed
\enddemo

\vskip0.4cm

Let $W\subseteq X$ be an irreducible component. 
Let $I$ be a semistable sheaf on $X$. It follows from Lemma 3 that there 
is a minimum subcurve $Z\subseteq X$ containing $W$ such that 
$\beta_I(W)=0$. We say that $I$ is \sl 
$W$-quasistable with respect to $E$ \rm if $\beta_I(Y)> 0$
for all proper subcurves $Y\subsetneqq X$ containing
$W$. It is clear that a semistable sheaf $I$ on $X$ 
is stable if and only if 
$I$ is $W$-quasistable with respect to $E$ 
for every irreducible component $W\subseteq X$. 

There is yet another interesting class of sheaves, this time 
independent of the 
choice of an irreducible component. 
We say that a semistable sheaf $I$ on 
$X$ is \sl quasistable with respect to 
$E$ \rm if there is an irreducible component $W\subseteq X$ 
such that $I$ is $W$-quasistable with respect to $E$. 
Note that a quasistable sheaf is
simple, an easy corollary of Prop. 1.

The notion of $W$-quasistability is not so easy to manage when dealing 
with families of curves. We shall often replace it with the equivalent, 
but more suitable, notion of $p$-stability. Let 
$p\in X$ be a non-singular point. We say that 
a semistable sheaf $I$ on $X$ is \sl 
$p$-quasistable with respect to $E$ \rm if $\beta_I(Y)>0$ 
for all proper subcurves $Y\subsetneqq X$ containing $p$. 

A semistable sheaf $I$ on $X$ is 
$X_i$-quasistable (for $i=1,\dots,n$) if and only if $I$ is 
$(r\underline{\delta}_i)$-quasistable (see \cite{\bf 13\rm , Sect. 4}), 
where $\underline{\delta}_i$ is the $n$-uple whose only non-zero 
component is the $i$-th component, with value 1.

We shall see in Subsect. 1.4 and Sect. 2 that the notion of 
quasistability is natural and useful.

\vskip0.4cm

\proclaim{Lemma 4} Let $Y,Z\subsetneqq X$ be proper subcurves 
covering $X$ such that $Y\wedge Z=\emptyset$, but $Y\cap Z\neq\emptyset$. 
Let $I$ (resp. $J$) be a torsion-free, rank 1 sheaf on $Y$ (resp. $Z$). 
Then, there is a non-split 
exact sequence of the form:
$$
0 @>>> J @>>> M @>>> I @>>> 0.
$$
\endproclaim

\demo{Proof} We need only show that 
$\text{Ext}^1_X(I,J)\neq 0$. Since $I$ and $J$ are 
torsion-free sheaves supported on subcurves $Y$ and $Z$, respectively, 
with $Y\wedge Z=\emptyset$, then 
$\underline{\text{Hom}}_X(I,J)=0$. Thus,
$$
\text{Ext}^1_X(I,J)=H^0(X,\underline{\text{Ext}}^1_X(I,J)).
$$
Of course, the topological support of $\underline{\text{Ext}}^1_X(I,J)$ 
is contained in $Y\cap Z$. Since $Y\cap Z\neq\emptyset$, there is 
$p\in Y\cap Z$. Let 
$\Cal O_p$ denote the local ring of $X$ at $p$, with maximal ideal 
$\Cal M_p$. Let $\Cal M_Y\subseteq \Cal O_p$ 
(resp. $\Cal M_Z\subseteq \Cal O_p$) 
be the ideal of $Y$ (resp. $Z$) at $p$. By 
hypothesis, $\Cal M_Y\cap \Cal M_Z=0$ and $\Cal M_Y+\Cal M_Z$ is a 
primary ideal of $\Cal M_p$. 
We need only 
show that $\text{Ext}^1_{\Cal O_p}(I_p,J_p)\neq 0$. Let
$$
(\frac{\Cal O_p}{\Cal M_Y})^{\oplus s_1} @>\phi >> 
(\frac{\Cal O_p}{\Cal M_Y})^{\oplus s_0} 
@>>> I_p @>>> 0\tag{4.1}
$$
be a presentation of $I_p$. Applying $\text{Ext}^1_{\Cal O_p}(-,J_p)$ to 
(4.1), we obtain a sequence,
$$
0 @>>> \text{Ext}^1_{\Cal O_p}(I_p,J_p) @>>> 
H^{\oplus s_0}  @>\phi^*\otimes\text{id}_H>> 
H^{\oplus s_1},\tag{4.2}
$$
where
$$
H:=\text{Ext}^1_{\Cal O_p}(\frac{\Cal O_p}{\Cal M_Y},J_p),
$$
and $\phi^*$ is the dual of $\phi$. The sequence (4.2) is 
exact, since $\text{Hom}_{\Cal O_p}(K,J_p)=0$ for all $\Cal O_p$-modules 
$K$ with $\Cal M_YK=0$. So, we need only show that 
$\phi^*\otimes\text{id}_H$ 
is not injective. Suppose, by contradiction, that 
$\phi^*\otimes\text{id}_H$ is injective. 
Since $H$ has finite length, then it follows from a standard argument that 
$\phi^*\otimes\text{id}_{\Cal O_p/\Cal M_p}$ is injective. Since $\phi^*$ 
is a homomorphism of free modules over the 
local ring $\Cal O_p/\Cal M_Y$, then 
$\phi^*$ is injective. 
It follows from (4.1) that
$$
\text{Hom}_{\Cal O_p}(I_p,\Cal O_p/\Cal M_Y)=0.
$$
Since $I$ is torsion-free, rank 1 on $Y$, then we have a contradiction. 
The proof is complete.\qed
\enddemo

\vskip0.4cm

The above lemma allows us to construct torsion-free, rank 1 sheaves on 
$X$ with ``prescribed'' Jordan-H\"older filtrations, as we shall see in 
Subsect. 1.3.

\vskip0.4cm

\parindent=0pt

\sl Example 5. \rm If $X$ is reducible, then there 
are semistable sheaves that are not simple. If $X$ has 
only two irreducible components, then every 
simple, semistable sheaf is quasistable. If $X$ has more than two 
components, then there might be simple, semistable sheaves that are 
not quasistable. For instance, suppose that 
there are connected subcurves $X_1, X_2, X_3\subseteq X$ covering $X$ 
such that $X_i\wedge X_j=\emptyset$ for $i\neq j$, but 
$X_1\cap X_2\neq\emptyset$ and $X_1\cap X_3\neq\emptyset$. Let 
$I_1,I_2,I_3$ be semistable, simple sheaves on $X_1,X_2,X_3$, 
respectively. By Lemma 4, since $X_1\cap X_2\neq\emptyset$ and 
$X_1\cap X_3\neq\emptyset$, then 
there is a non-split exact sequence of the form:
$$
0 @>>> I_2\oplus I_3 @>>> I @>>> I_1 @>>> 0,
$$
whose push-out to $I_2$ (resp. $I_3$) is 
a non-split exact sequence of the form:
$$
0 @>>> I_2 @>>> I_{X_1\cup X_2} @>>> I_1 @>>> 0 \text{\  \  (resp. \  } 
0 @>>> I_3 @>>> I_{X_1\cup X_3} @>>> I_1 @>>> 0\text{)}. 
$$
It is clear that $I$, $I_{X_1}$, $I_{X_1\cup X_2}$ and $I_{X_1\cup X_3}$ are 
semistable. So $I$ is not quasistable. Moreover, since the latter 
exact sequences are non-split, then $I_{X_1\cup X_2}$ and $I_{X_1\cup X_3}$ 
are simple by statement (2) of Lemma 2. It follows now from 
statement (1) of Lemma 2 that $I$ is simple. We have thus produced a 
simple, semistable sheaf $I$ that is not quasistable.

\vskip0.4cm

\parindent=0pt

\subheading{1.3. Jordan-H\"older filtrations} 
Let $I$ be a semistable sheaf on $X$ with respect to $E$. 
We describe now a filtration of $I$. To start with, we let $I_0:=I$ 
and $Z_0:=X$. 
Let $Y_0\subseteq X$ be a non-empty subcurve such that $I_{Y_0}$
is stable with respect to $\left. E\right|_{Y_0}$. 
Let $I_1:=\text{ker}(I @>>> I_{Y_0})$. Clearly, the
sheaf $I_1$ is torsion-free, rank 1 on $Z_1:=Y^c_0$, 
and semistable with respect to $\left.E\right|_{Z_1}$, if not zero. 
Repeating the above procedure with $I_1$, in the place of
$I$, and so on, we end up with filtrations:
$$
\left\{\aligned 
0=& I_{q+1}\subsetneqq I_q \subsetneqq \dots \subsetneqq I_1\subsetneqq 
I_0=I,\\
\emptyset =& Z_{q+1}\subsetneqq Z_q\subsetneqq\dots\subsetneqq 
Z_1\subsetneqq Z_0=X,
\endaligned\right.
$$
with the following properties:
\roster
\item for $i=0,\dots,q$, the sheaf 
$I_i$ is torsion-free, rank 1 on the subcurve $Z_i\subseteq X$, 
and is semistable
with respect to $\left.E\right|_{Z_i}$;
\item for $i=0,\dots,q$, the quotient 
$I_i/I_{i+1}$ is torsion-free, rank 1 on the subcurve $Y_i:=Z_i-Z_{i+1}$, 
and is stable with respect to $\left.E\right|_{Y_i}$.
\endroster

\parindent=12pt

We call the above filtration of $I$ a 
\sl Jordan-H\"older filtration\rm . 
A Jordan-H\"older filtration depends on the choices
made in its construction, but
$$
\text{Gr}(I):=I_0/I_1\oplus I_1/I_2\oplus\dots\oplus I_q/I_{q+1}
$$
depends only on $I$, by the Jordan-H\"older theorem. 
In particular, the collection of
subcurves $\{Y_0,\dots,Y_q\}$ covering $X$ depends only on
$I$. It is clear that $\text{Gr}(I)$ is torsion-free, rank 1, and 
semistable with respect to $E$. In addition, we have that 
$\text{Gr}(\text{Gr}(I))=\text{Gr}(I)$. If $I$ is stable, 
then $\text{Gr}(I)=I$. We shall say that two semistable sheaves, $I$ and $J$, 
are \sl 
Jordan-H\"older equivalent \rm (\sl JH-equivalent\rm , for short) 
if $\text{Gr}(I)\cong \text{Gr}(J)$. 

\vskip0.4cm

\proclaim{Proposition 6} Let $I$ be a semistable sheaf on $X$. Let
$$
\left\{\aligned 
0=& I_{q+1} \subsetneqq I_q 
\subsetneqq \dots \subsetneqq I_1 \subsetneqq I_0=I,\\
\emptyset =& Z_{q+1} \subsetneqq Z_q \subsetneqq \dots 
\subsetneqq Z_1 \subsetneqq Z_0=X
\endaligned\right.
$$
be a Jordan-H\"older filtration of $I$. Let $W\subseteq X$ be an 
irreducible component. Then, the following statements are equivalent:
\roster 
\item $I$ is $W$-quasistable;
\item $I_i$ is $W$-quasistable for $i=0,\dots,q$;
\item $W\subseteq Z_q$, and the short exact sequence:
$$
0 @>>> I_{i+1} @>>> I_i @>>> \frac{I_i}{I_{i+1}} @>>> 0
$$
is not split for $i=0,\dots,q-1$.
\endroster
\endproclaim

\demo{Proof} See \cite{\bf 13\rm , Prop. 5}.\qed
\enddemo

\vskip0.4cm

The following theorem shows that in every Jordan-H\"older 
equivalence class there is a 
$W$-quasistable sheaf, for any fixed irreducible component $W\subseteq X$.

\vskip0.4cm

\proclaim{Theorem 7} Assume that $X$ is connected. 
Let $Y_0,\dots,Y_q\subseteq X$ be subcurves covering $X$ such that 
$Y_i\wedge Y_j=\emptyset$ for $i\neq j$. Let $J_0,\dots,J_q$ be stable sheaves 
on $Y_0,\dots, Y_q$, respectively. Let 
$W\subseteq X$ be an irreducible component. Then, there is a 
$W$-quasistable sheaf $I$ on $X$ such that 
$\text{\rm Gr}(I)\cong J_0\oplus\dots\oplus J_q$.
\endproclaim

\demo{Proof} We may assume that $W\subseteq Y_q$. Since $X$ is connected, 
we may also assume that $(Y_q\cup\dots\cup Y_i)\cap Y_{i-1}\neq\emptyset$ 
for $i=1,\dots,q$. Let $Z_i:=Y_i\cup\dots\cup Y_q$ for $i=0,\dots,q$. Let 
$I_q:=J_q$. We construct recursively 
a torsion-free, rank 1 sheaf 
$I_i$ on $Z_i$ for $i=q-1,\dots,0$ as follows: 
suppose that we are given $I_q,\dots,I_i$ for a certain 
$i\in\{1,\dots,q\}$; then, 
let $I_{i-1}$ be the middle sheaf in a non-split exact sequence of the 
form:
$$
0 @>>> I_i @>>> I_{i-1} @>>> J_{i-1} @>>> 0,
$$
whose existence is guaranteed by Lemma 4. It is clear that 
we have a Jordan-H\"older filtration
$$
\left\{\aligned 
0 & =:I_{q+1} \subsetneqq I_q \subsetneqq \dots \subsetneqq I_1 
\subsetneqq I_0=:I,\\
\emptyset & =:Z_{q+1} \subsetneqq Z_q \subsetneqq \dots 
\subsetneqq Z_1 \subsetneqq Z_0=X
\endaligned\right.
$$
of $I:=I_0$. Of course, $\text{Gr}(I)=J_0\oplus\dots\oplus J_q$. In addition, 
it follows from Prop. 6 that $I$ is $W$-quasistable. 
The proof is complete.\qed
\enddemo

\vskip0.4cm

\parindent=0pt

\subheading{1.4. The cohomological characterizations} 
Assume that $k$ is algebraically closed. Let $\omega$ be the dualizing
sheaf on $X$. Recall that $\omega$ is simple, torsion-free, rank 1. 

\parindent=12pt

\vskip0.4cm

For the next three results, let $E$ be a vector bundle on $X$ with rank 
$r>0$ and degree $-rd$, for some integer $d$.

\vskip0.4cm

\proclaim{Theorem 8} Let $I$ be a torsion-free, rank 1
sheaf on $X$ with $\chi(I)=d$. Then, $I$ is semistable with respect to $E$
if and only if there is a vector bundle $F$ on $X$, with rank $mr$ and 
$\det F\cong(\det E)^{\otimes m}$ for some $m>0$, such that
$$
h^0(X,I\otimes F)=h^1(X,I\otimes F)=0.
$$
\endproclaim

\demo{Proof} As in \cite{\bf 28\rm , Lemma 3.1, p. 166} for 
the ``only if'' part, and \cite{\bf 28\rm , Lemma 8.3, p. 195} 
for the ``if'' part.\qed
\enddemo

\vskip0.4cm

\proclaim{Theorem 9} Let $p\in X$ be a non-singular point. 
Let $I$ be a semistable sheaf on $X$ with respect to
$E$. Then, $I$ is $p$-quasistable if and only if there
is a vector bundle $F$ on $X$, with rank $mr$ and 
$\det F\cong(\det E)^{\otimes m}\otimes\Cal O_X(-p)$ for some 
$m>0$, such that:
\roster
\item $h^0(X,I\otimes F)=0$ and $h^1(X,I\otimes F)=1$;
\item the unique (modulo $k^*$) non-zero homomorphism 
$I @>>> F^*\otimes\omega$ is injective.
\endroster 
\endproclaim

\demo{Proof} The ``if'' part follows directly from 
\cite{\bf 13\rm , Thm. 6}. The ``only if'' part follows from loc. cit. and 
a standard argument.\qed
\enddemo

\vskip0.4cm

\proclaim{Corollary 10} Let $I$ be a semistable sheaf on $X$ with respect 
to $E$. Then, $I$ is quasistable (resp. stable) if and only if 
for some (resp. for every) non-singular point $p\in X$ there is a vector 
bundle $F$ on $X$, with rank $mr$ and $\det F\cong(\det E)^{\otimes m}
\otimes\Cal O_X(-p)$ for some $m>0$, such that:
\roster 
\item $h^0(X,I\otimes F)=0$ and $h^1(X,I\otimes F)=1$;
\item the unique (modulo $k^*$) non-zero homomorphism $I @>>> 
F^*\otimes\omega$ is injective.
\endroster
\endproclaim

\demo{Proof} By definition, $I$ is
quasistable (resp. stable) 
if and only if $I$ is $p$-quasistable for some (resp. for every) non-singular 
point $p\in X$. Thus, the 
corollary follows immediately from Thm. 9.\qed
\enddemo

\vskip0.4cm

For every coherent sheaf $I$ on $X$, we let 
$$
e(I):=\max_{p\in X} \dim_k I(p).
$$

\vskip0.4cm

\proclaim{Proposition 11} Let $W\subseteq X$ be an irreducible 
component. 
If $I$ is a simple, torsion-free, rank 1 sheaf on $X$, then there is a 
vector bundle $E$ on $X$ of rank 
$r:=\max(e(\underline{\text{\rm Hom}}_X(I,\omega)),2)$ 
such that $I$ is $W$-quasistable 
with respect to $E$.
\endproclaim

\demo{Proof} Applying the proof of \cite{\bf 12\rm , Prop. 1}, we get a 
vector bundle $F$ on $X$ with 
rank $r$ and degree $-r\chi(I)-1$, such that
$$
h^0(X,I\otimes F)=0\text{\  \  and \  \  } h^0(X,I\otimes F)=1,
$$
and the unique (modulo $k^*$) non-zero homomorphism 
$\lambda\: I @>>> F^*\otimes\omega$ is injective. 
(Note: we use that $I$ is simple 
in order to apply the 
proof of loc. cit..) Let $p\in W$ be a non-singular point of $X$. Let
$$
E:=\text{ker}(F^* @>>> F^*(p) @>q>> k)^*,
$$
where $q$ is a $k$-linear surjective homomorphism 
such that $q\circ\lambda(p)\neq 0$. 
(Note that we chose implicitly a trivialization of $\omega$ at $p$ to 
consider the latter composition.) It is clear 
that $E$ has rank $r$ and degree $-r\chi(I)$. In fact, 
$\det E\cong(\det F)\otimes\Cal O_X(p)$. In addition, by our choice of 
$q$, we have that
$$
h^0(X,I\otimes E)=h^1(X,I\otimes E)=0.
$$
By Thm. 8, the sheaf $I$ is semistable with respect to $E$. Moreover, 
it follows from Thm. 9 that $I$ is $W$-quasistable with respect to $E$. 
The proof is complete.\qed
\enddemo

\vskip0.4cm

\parindent=0pt

\sl Remark 12. \rm Given an irreducible component $W\subseteq X$, it 
seems that a 
simple, torsion-free, rank 1 sheaf $I$ on $X$ is $W$-quasistable with 
respect to a line bundle. 
At least, I can show the latter statement in the following cases:
\roster
\item $X$ is Gorenstein and $I$ is invertible;
\item $X$ has only ordinary nodes as singularities;
\item $X$ has at most two irreducible components.
\endroster 
Since this statement is not central to the present work, I 
will omit the proof of 
it in the above three cases. 

\vskip0.4cm

\parindent=12pt

\proclaim{Corollary 13} If $I$ is a simple, torsion-free, rank 1 
sheaf on $X$, then there is a vector bundle $E$ on $X$ of 
rank $nr$, where $r:=\max(e(\underline{\text{\rm Hom}}_X(I,\omega)),2)$, 
and $n$ is the number of irreducible components of $X$, 
such that $I$ is stable with respect to $E$.
\endproclaim

\demo{Proof} By Prop. 11, for each $i=1,\dots,n$ there is a 
vector bundle $E_i$ on $X$ of rank $r$ such that $I$ is $X_i$-quasistable 
with respect to $E_i$. Then, $I$ is stable with respect to 
the direct sum, $E:=E_1\oplus\dots\oplus E_n$.\qed
\enddemo

\vskip0.4cm

\parindent=0pt

\sl Remark 14. \rm If the expectation stated in 
Rmk. 12 is confirmed, then it will follow that every simple, 
torsion-free, rank 1 sheaf on $X$ is 
stable with respect to a rank $n$ vector bundle. The expected upper bound 
$n$ on the rank of the polarization is sharp: we can easily construct 
an example of a curve $X$ with $n$ irreducible components, and 
a simple, torsion-free, rank 1 sheaf $I$ on $X$ 
such that $I$ is not stable with respect to any 
vector bundle of rank less than $n$. 

\vskip0.4cm

\parindent=0pt

\sl Remark 15. \rm The cohomological characterizations, 
Thm. 8 and Thm. 9, remain valid if the 
curve $X$ is defined over an infinite field $k$. 
(In Thm. 9 we have to assume that $p$ is $k$-rational.) 
If $X$ is defined over a finite field $k$, 
then the said theorems remain valid after a finite field 
extension of $k$. 
To see that, we proceed as follows. Let $\bar k$ denote the algebraic 
closure of $k$. First, we use Thm. 8 (resp. Thm. 9) 
to get a vector bundle $G$ on $\bar X:=X\times_k \bar k$, with rank $mr$ and 
$\det G\cong L\otimes_k \bar k$, with the required properties, where 
$L:=(\det E)^{\otimes m}$ 
(resp. $L:=(\det E)^{\otimes m}\otimes\Cal O_X(-p)$). Second, we observe 
that $G$ fits in the middle of a short exact sequence of the form:
$$
0 @>>> (A^{\otimes -c})^{\oplus (mr-1)}\otimes_k \bar k @>>> G @>>> 
L\otimes A^{\otimes c(mr-1)}\otimes_k \bar k @>>> 0,\tag{15.1}
$$
where $A$ is an ample invertible sheaf on $X$, and $c>>0$. Let
$$
V:=H^1(X,(A^{\otimes -cmr})^{\oplus (mr-1)}\otimes L^*).
$$
The affine space $\text{\bf A}(V)$ parametrizes exact sequences of the 
form (15.1). Thus, $G$ is ``represented'' by a geometric point 
$g\in \text{\bf A}(V)$. 
Since the properties required of $G$ in Thm. 8 (resp. Thm. 9) are ``open'', 
then there is an open neighbourhood 
$U\subseteq\text{\bf A}(V)$ of $g$ 
where these properties remain valid. If $k$ is infinite, then there is 
a $k$-rational point $f\in U$. 
If $k$ is finite, after replacing $k$ by a finite field extension of 
sufficiently high degree, if necessary, there will be a 
$k$-rational point $f\in U$. The point $f$ 
corresponds to an exact sequence of the form:
$$
0 @>>> (A^{\otimes -c})^{\oplus (mr-1)} @>>> F @>>> 
L\otimes A^{\otimes c(mr-1)} @>>> 0,
$$
where the vector bundle $F$ satisfies the requirements of Thm. 8 (resp. 
Thm. 9).

\parindent=0pt

\vskip0.4cm

\sl Remark 16. \rm As it could be expected, we can find the vector bundle 
$F$ in Thm. 8 and Thm. 9 with rank restricted to a certain range, depending 
only on numerical invariants. The existence of such range follows 
easily from the fact that the family of all 
semistable sheaves on $X$, with respect to $E$, is bounded. 

\parindent=0pt

\vskip0.4cm

\subheading{1.5. Families} Let $f\: X @>>> S$ be
a flat, projective morphism whose geometric fibres are curves. 
Let $\Cal I$ be an $S$-flat coherent sheaf on $X$. We say that 
$\Cal I$ is \sl relatively torsion-free \rm  (resp.
\sl rank 1, \rm  resp. \sl simple\rm ) over $S$ if $\Cal I(s)$ is
torsion-free (resp. rank 1, resp. simple) for every 
geometric point $s\in S$.

\parindent=12pt

Let $\Cal F$ be an $S$-flat coherent sheaf (resp. a vector
bundle) on $X$. 
We say that $\Cal F$ has \sl relative Euler characteristic $d$ \rm
(resp. \sl relative degree $d$\rm) over $S$ if $\Cal
F(s)$ has Euler characteristic $d$ (resp. degree $d$) 
for every geometric point $s\in S$.

Let $\Cal E$ be a vector bundle on $X$
of rank $r>0$ and relative degree $-rd$ over $S$, for a certain integer 
$d$. We call such an $\Cal E$ a \sl 
relative polarization \rm on $X$ over $S$. A relatively torsion-free, rank 1
sheaf $\Cal I$ on $X$ over $S$ is \sl relatively stable \rm  (resp. \sl
semistable\rm, resp. \sl quasistable\rm) \sl with respect
to $\Cal E$ \rm over $S$ if $\Cal I(s)$ is
stable (resp. semistable, resp. quasistable) with respect to $\Cal E(s)$ for
every geometric point $s\in S$. Let $\sigma\: S @>>> X$ be a section of $f$ 
through the $S$-smooth locus of $X$. A relatively torsion-free, rank 1 
sheaf $\Cal I$ on $X$ over $S$ is \sl relatively $\sigma$-quasistable 
with respect to $\Cal E$ \rm over $S$ if $\Cal I(s)$ is 
$\sigma(s)$-quasistable with respect to $\Cal E(s)$ for every 
geometric point $s\in S$.

\vskip0.8cm

\parindent=12pt

\heading
2. The valuative criteria
\endheading

\parindent=12pt

Let $S$ be the spectrum of a discrete valuation ring
$R$. Let $s$ be the special point of $S$, and $\eta$ be
its generic point. Let $f\: X @>>> S$ be a flat, projective 
morphism whose geometric fibres are curves. Assume that the irreducible
components of the special fibre $X(s)$ are geometrically integral. 

If $\Cal I$ is a relatively torsion-free, rank 1
sheaf on $X$ over $S$, and $Y\subseteq X(s)$ is a subcurve, 
then we denote by $\Cal I^Y$ the
kernel of the canonical surjective homomorphism,
$$
\Cal I \twoheadrightarrow \Cal I(s)_Y.
$$
It is clear that the inclusion map, $\lambda\:\Cal I^Y
\hookrightarrow \Cal I$, is
an isomorphism on $X\setminus Y$. In addition, it can be shown, using an 
argument analogous to the one found in 
\cite{\bf 21\rm , Prop. 6, p. 100}, that 
$\Cal I^Y$ is relatively torsion-free, rank 1 
on $X$ over $S$. 

The following lemma is the main technical tool of this section.

\vskip0.4cm

\proclaim{Lemma 17} Let $Y\subseteq X(s)$ be a subcurve. Let
$$
\dots\subseteq\Cal I^i\subseteq\dots\subseteq\Cal I^1
\subseteq\Cal I^0:=\Cal I
$$
be an infinite filtration of $\Cal I$ with quotients
$$
\frac{\Cal I^i}{\Cal I^{i+1}}=\Cal I^i(s)_Y
$$
for every $i\geq 0$. If $R$ is complete, and 
$\Cal I^i(s)$ decomposes at $Y$ for every $i\geq 0$, 
then there is an $S$-flat quotient $\Cal F$ of
$\Cal I$ on $X$ such that $\Cal F(s)=\Cal I(s)_Y$.
\endproclaim

\demo{Proof} The proof is very similar to the proof of \cite{\bf 
21\rm , Lemma 2, p. 106}, hence we just outline the construction of 
$\Cal F$. Let $\pi$ be a generator of the maximal ideal of $R$. 
For every $i\geq 1$, let $S_i:=\text{Spec }(R/(\pi^{i-1}))$ 
and $X_i:=X\times_S S_i$. Let
$$
F_i:=\text{coker}(\left. \Cal I^i\right|_{X_i} @>>> \left. \Cal I\right|_
{X_i}).
$$
We can show the 
following two properties:
\roster 
\item $\left. F_i\right|_{X_j}=F_j$ if $i\geq j\geq 1$;
\item $F_i$ is $S_i$-flat for every $i\geq 1$.
\endroster 
By Grothendieck's existence theorem \cite{\bf 16\rm , III-1-5.1.7},
since $R$ is complete, there is a quotient $\Cal F$ of
$\Cal I$ on $X$ such that $\Cal F$ is the inverse
limit of the $F_i$. Since each $F_i$ is $S_i$-flat, then
$\Cal F$ is $S$-flat. Moreover, $\Cal F(s)=F_1=\Cal I(s)_Y$. The 
outline is complete.\qed
\enddemo

\vskip0.4cm

Fix an integer $d$, and 
let $\Cal E$ be a vector bundle on $X$ of
rank $r>0$ and relative degree $-rd$ over $S$. We fix $\Cal E$ as 
our relative polarization (see Subsect. 1.5).

Let $Z\subseteq X(s)$ be a subcurve. Let $\Cal I$ 
be a relatively torsion-free, rank 1 sheaf on $X$ over $S$. Let 
$$
\Cal I':=\text{ker}(\Cal I @>>> \Cal I(s)_Z) \text{\  \  and \  \  } 
\Cal I'':=\text{ker}(\Cal I' @>>> \Cal I'(s)_{Z^c}).
$$
We claim that $\Cal I''=\pi\Cal I$. 
In fact, the claim is clearly true 
on the open subset $U\subseteq X$ 
obtained by excluding the singular 
points of $X(s)$. But $U$ contains all
points of depth less than 2 in $X$. Thus, the 
canonical homomorphism $\Cal K @>>> \left.\Cal K\right|_U$ is an
isomorphism for every relatively torsion-free, rank 1 sheaf $\Cal K$
on $X$ over $S$. 
The claim follows easily now. Since $\Cal I''=\pi \Cal I$, then 
$\Cal I''\cong \Cal I$. Let $\lambda\: \Cal I' \hookrightarrow 
\Cal I$ denote the 
inclusion homomorphism, and let $\mu\: \Cal I \hookrightarrow \Cal I'$ 
denote the embedding induced by the inclusion $\Cal I''\subseteq\Cal I'$ and 
the isomorphism $\Cal I''\cong\Cal I$. It is clear that 
$\lambda(s)$ and $\mu(s)$ induce the exact sequences:
$$
\aligned 
0 & @>>> \Cal I'(s)_{Z^c} @>>> \Cal I(s) @>>> \Cal I(s)_Z @>>> 0,\\
0 & @>>> \Cal I(s)_Z @>>> \Cal I'(s) @>>> \Cal I'(s)_{Z^c} @>>> 0,
\endaligned\tag{17.1}
$$
respectively.

\vskip0.4cm

\proclaim{Lemma 18} Let $\Cal I$ be a relatively semistable sheaf 
on $X$ over $S$. If $Z\subseteq X(s)$ is a subcurve, then $\Cal I^Z$ 
is relatively 
semistable on $X$ over $S$ if and only if 
$\Cal I(s)_Z$ is semistable with respect to $\left.\Cal E(s)\right|_Z$. If 
that is the case, then $\Cal I^Z(s)$ is JH-equivalent to $\Cal I(s)$.
\endproclaim

\demo{Proof} Of course, $\Cal I^Z$ is relatively semistable if and only 
if $\Cal I^Z(s)$ is semistable. The lemma follows easily now from 
considering the two exact sequences in (17.1). 
The proof is complete.\qed
\enddemo

\vskip0.4cm

\proclaim{Proposition 19} Let $\Cal I$, $\Cal J$ be
relatively semistable sheaves on $X$ over $S$ such that 
$\Cal I(\eta)\cong\Cal J(\eta)$. Then, 
$\Cal I(s)$ and $\Cal J(s)$ are JH-equivalent. In particular, if $\Cal I(s)$ 
is stable, then $\Cal I(s)\cong\Cal J(s)$.
\endproclaim

\demo{Proof} By the same argument used in \cite{\bf 4\rm , Lemma 7.8, 
p. 100}, there is a homomorphism 
$\lambda\: \Cal I @>>> \Cal J$ such that $\lambda(\eta)$
is an isomorphism, and $\lambda(s)$ is non-zero. Let 
$Y\subseteq X(s)$ be the subcurve such that 
$\Cal I(s)_Y\cong\text{im}(\lambda(s))$. 
Since $\lambda(s)$ is non-zero, then $Y$ is non-empty. 
If $Y=X(s)$, then $\lambda(s)$ is an embedding. In this case, 
since $\lambda(\eta)$ is an 
isomorphism, then so is $\lambda(s)$. We may thus assume that $Y$ is a 
proper subcurve of $X$. Since $\Cal I(s)$ and $\Cal J(s)$ are semistable, 
then so is $\Cal I(s)_Y$. Since $\Cal I(s)_Y$ and $\Cal J(s)$ are 
semistable, then so is $\Cal J(s)_Z$, where $Z:=Y^c$. 
Let $\Cal J_1:=\text{ker}(\Cal J @>>> \Cal J(s)_Z)$. Since 
$\Cal J(s)_Z$ is semistable, then $\Cal J_1$ is relatively semistable on 
$X$ over $S$, and $\Cal J_1(s)$ is JH-equivalent to 
$\Cal J(s)$, by Lemma 18. Moreover, 
$\lambda$ factors through $\Cal J_1$ by construction. Applying the same 
procedure described above to the induced $\lambda_1\: \Cal I @>>> \Cal J_1$, 
in the place of $\lambda$, and so on, we will eventually find a 
relatively semistable subsheaf $\Cal J_m\subseteq \Cal J$ 
such that $\Cal J_m(s)$ is JH-equivalent to $\Cal J(s)$, and the induced 
homomorphism, $\lambda_m\: \Cal I @>>> \Cal J_m$, 
is an isomorphism. The proof is complete.\qed
\enddemo

\vskip0.4cm

\proclaim{Proposition 20} Let $\sigma\: S @>>> X$ be a
section of $f$ through the $S$-smooth locus of $X$. Let $\Cal I$
and $\Cal J$ be relatively $\sigma$-quasistable sheaves on
$X$ over $S$ with respect to $\Cal E$. If $\Cal I(\eta)\cong\Cal
J(\eta)$, then $\Cal I\cong\Cal J$.
\endproclaim

\demo{Proof} As in the proof of \cite{\bf 4\rm , Lemma 7.8, p. 100}, 
there are homomorphisms $\lambda\: \Cal I @>>> \Cal J$ and
$\mu\: \Cal J @>>> \Cal I$ such that $\lambda(\eta)$ and
$\mu(\eta)$ are isomorphisms, and $\lambda(s)$ and
$\mu(s)$ are non-zero. Let $Y,Z\subseteq X(s)$ be subcurves 
such that $\Cal I(s)_Y\cong\text{im}(\lambda(s))$ and 
$\Cal J(s)_Z\cong\text{im}(\mu(s))$. Since
$\lambda(s)$ and $\mu(s)$ are non-zero, then $Y$ and $Z$
are non-empty. Since $\Cal I(s)$
and $\Cal J(s)$ are semistable, then so are $\Cal
I(s)_{Z^c}$ and $\Cal J(s)_{Y^c}$. Since $\Cal I(s)$ and
$\Cal J(s)$ are $\sigma(s)$-quasistable, then
$\sigma(s)\notin Y^c\cup Z^c$. 
So, $\sigma(s)\in Y\wedge Z$. Since $Y\wedge Z\neq \emptyset$, then the
composition $\mu(s)\circ\lambda(s)$ is not zero.
Since $\Cal I(s)$ is simple, then $\mu(s)\circ\lambda(s)$ is a 
homothety. Thus, $\lambda(s)$ is
injective and, since $\lambda(\eta)$ is an isomorphism, then 
$\lambda$ is an isomorphism. The proof
is complete.\qed
\enddemo

\vskip0.4cm

\proclaim{Lemma 21} Let $\Cal I$ be a relatively torsion-free, rank 1 
sheaf on $X$ over $S$. Let $Z\subseteq X(s)$ be a subcurve such that 
$\Cal I(s)$ decomposes at $Z$. Then, for any subcurve
$Y\subseteq X$, if $\Cal I^Z(s)$ decomposes at $Y$, then so does 
$\Cal I(s)$.
\endproclaim 

\demo{Proof} The lemma follows easily from considering the two exact 
sequences in (17.1). The proof is complete.\qed
\enddemo

\vskip0.4cm

Let $I$ be a torsion-free, rank 1 sheaf on $X(s)$. 
If $Y,Z\subseteq X(s)$ are subcurves with $Y\wedge Z=\emptyset$, we let
$$
\delta_I(Y,Z):=\chi(I_Y)+\chi(I_Z)-\chi(I_{Y\cup Z}).
$$
If $Z'\subseteq Z$ is a subcurve, then
$\delta_I(Y,Z')\leq\delta_I(Y,Z)$ by Lemma 3. 
In particular, $\delta_I(Y,Z)\geq 0$, with equality if and only if $I_{Y\cup
Z}=I_Y\oplus I_Z$.

\vskip0.4cm

\proclaim{Lemma 22} Let $\Cal I$ be a relatively torsion-free, rank 1 
sheaf on $X$ with relative Euler characteristic $d$ 
over $S$. Let $Y,Z\subseteq X$ be subcurves. Then,
$$
\beta_{\Cal I^Z(s)}(Y)+\beta_{\Cal I(s)}(Z)\geq \beta_{\Cal I(s)}(Y\wedge Z)
+\beta_{\Cal I(s)}(Y\cup Z),
$$
with equality if and only if $\delta_{\Cal I^Z(s)}(Y\wedge Z,Z^c)=
\delta_{\Cal I^Z(s)}(Y\wedge Z,Y\wedge Z^c)$.
\endproclaim

\demo{Proof} The lemma follows from considering 
the exact sequences in (17.1), together with 
an argument analogous to the 
one found in \cite{\bf 21\rm , Lemma 1, p. 105}. 
The proof is complete.\qed
\enddemo

\vskip0.4cm

Let $I$ be a torsion-free, rank 1 sheaf on $X(s)$ with 
$\chi(I)=d$. It follows 
easily from Lemma 3 that there is a maximum subcurve $Z\subseteq X(s)$ among 
the subcurves $W\subseteq X(s)$ with minimum $\beta_I(W)$. It is 
clear that $I$ is semistable if and only if $Z=X(s)$.

\vskip0.4cm

\proclaim{Lemma 23} Let $\Cal I$ be a relatively torsion-free, rank 1 
sheaf on $X$ of relative Euler characteristic $d$ over $S$. Let 
$Z\subseteq X(s)$ be the maximum subcurve among the subcurves 
$W\subseteq X(s)$ with minimum $\beta_{\Cal I(s)}(W)$. Then, 
$\beta_{\Cal I^Z(s)}(Y)\geq\beta_{\Cal I(s)}(Z)$ for every 
subcurve $Y\subseteq X(s)$, with equality only if $Y\subseteq Z$. 
Moreover, we have that $\beta_{\Cal I^Z(s)}(Z)=\beta_{\Cal I(s)}(Z)$ if
and only if $\Cal I^Z(s)$ decomposes at $Z$.
\endproclaim

\demo{Proof} Let $Y\subseteq X(s)$ be a subcurve. 
Since $\beta_{\Cal I(s)}(Z)$ is minimum, it follows from 
Lemma 22 that
$$
\beta_{\Cal I^Z(s)}(Y)\geq\beta_{\Cal I(s)}(Y\wedge Z)-\beta_{\Cal I(s)}
(Z)+\beta_{\Cal I(s)}(Y\cup Z)\geq\beta_{\Cal I(s)}(Z),
$$
with equality if and only if
$$
\aligned
\delta_{\Cal I^Z(s)}(Y\wedge Z,Z^c)=&\delta_{\Cal I^Z(s)}(Y\wedge Z,Y\wedge 
Z^c),\\
\beta_{\Cal I(s)}(Y\cup Z)=&\beta_{\Cal I(s)}(Z),\\
\beta_{\Cal I(s)}(Y\wedge Z)=&\beta_{\Cal I(s)}(Z).
\endaligned
$$
Since $Z$ is maximum among the subcurves $W\subseteq X(s)$ with 
minimum $\beta_{\Cal I(s)}(W)$, then the above 
middle equality occurs if and 
only if $Y\subseteq Z$. The first statement of the
lemma is proved. If $Y=Z$, then the above last two equalities are 
obviously satisfied, whereas the first equality is
satisfied if and only if $\delta_{\Cal I^Z(s)}(Z,Z^c)=0$. 
The proof is complete.\qed
\enddemo

\vskip0.4cm

\proclaim{Lemma 24} Let $\Cal I$ be a relatively 
semistable sheaf on $X$ over $S$ with 
respect to $\Cal E$. Let $W\subseteq X(s)$ 
be an irreducible component, and 
let $Z\subseteq X(s)$ be the minimum subcurve containing $W$ such that 
$\beta_{\Cal I(s)}(Z)=0$. Then, $\Cal I^Z$ is also relatively 
semistable on $X$ over $S$ with respect to $\Cal E$. Moreover, 
if $Z'\subseteq X(s)$ is the minimum subcurve containing $W$ such that 
$\beta_{\Cal I^Z(s)}(Z')=0$, then $Z'\supseteq Z$, with equality 
$Z'=Z$ if and only if $\Cal I^Z(s)$ decomposes at $Z$.
\endproclaim

\demo{Proof} The first statement is a direct application of Lemma 18. 
As for the second statement, since $\Cal I(s)$ is semistable, 
$\beta_{\Cal I(s)}(Z)=0$ and 
$\beta_{\Cal I^Z(s)}(Z')=0$, 
it follows from Lemma 22 that
$$
0=\beta_{\Cal I^Z(s)}(Z')\geq 
\beta_{\Cal I(s)}(Z'\wedge Z)+\beta_{\Cal I(s)}(Z'\cup Z)\geq 0.
$$
So, $\beta_{\Cal I(s)}(Z'\wedge Z)=0$ and, by Lemma 22 again,
$$
\delta_{\Cal I^Z(s)}(Z'\wedge Z,Z^c)=\delta_{\Cal I^Z(s)}(Z'\wedge Z,
Z'\wedge Z^c).\tag{24.1}
$$
Since $Z'\supseteq W$, and $Z$ is the minimum subcurve containing $W$ with 
$\beta_{\Cal I(s)}(Z)=0$, then $Z'\supseteq Z$. 
The rest of the 
second statement follows now from (24.1).\qed
\enddemo

\vskip0.4cm

Let $I_{\eta}$ be a torsion-free, rank 1 sheaf on $X(\eta)$. We say that a 
relatively torsion-free, rank 1 sheaf $\Cal I$ on $X$ over $S$ such that 
$\Cal I(\eta)\cong I_{\eta}$ is an \sl extension of $I_{\eta}$\rm .

\vskip0.4cm

\proclaim{Theorem 25} Let $I_{\eta}$ be a torsion-free, rank 1 sheaf 
on $X(\eta)$. Then, the following statements hold:
\roster
\item There is an extension $\Cal I$ of $I_{\eta}$.
\item If $I_{\eta}$ is simple, then there is an extension $\Cal I$ of 
$I_{\eta}$ that is relatively simple over $S$.
\item If $I_{\eta}$ is (simple and) semistable 
with respect to $\Cal E(\eta)$, then 
there is an extension $\Cal I$ of $I_{\eta}$ that is 
relatively (simple and) semistable over $S$ with respect to $\Cal E$.
\item Let $\sigma\: S @>>> X$ be a section of $f$ through the $S$-smooth 
locus of $X$. If $I_{\eta}$ is $\sigma(\eta)$-quasistable with respect to 
$\Cal E(\eta)$, then there is an extension $\Cal I$ of $I_{\eta}$ that is 
relatively $\sigma$-quasistable over $S$ with respect to $\Cal E$.
\item If $I_{\eta}$ is quasistable with respect to $\Cal E(\eta)$, then 
there is an extension $\Cal I$ of $I_{\eta}$ that is relatively 
quasistable over $S$ with respect to $\Cal E$.
\endroster
\endproclaim

\demo{Proof} Statement (1) follows immediately from the same argument used in 
\cite{\bf 4\rm , Lemma 7.8 (i), p. 100}.

We prove (2) now. By (1), we may pick 
an extension $\Cal I$ of $I_{\eta}$. 
If $\Cal I(s)$ is simple, then we are done. 
If not, it follows from Prop. 1 that 
there is a non-empty, proper subcurve 
$Z\subsetneqq X(s)$ 
such that $\Cal I(s)$ decomposes at $Z$. In this case, let 
$$
\Cal I^1:=\text{ker}(\Cal I @>>> \Cal I(s)_Z) \text{\  \  and \  \  } 
\Cal I^{-1}:=\text{ker}(\Cal I @>>> \Cal I(s)_{Z^c}).
$$
By Lemma 21, the set $\Cal C^1$ (resp. $\Cal C^{-1}$) 
of subcurves $Y\subseteq X(s)$ such that 
$\Cal I^1(s)$ (resp. $\Cal I^{-1}(s)$) decomposes at $Y$ is contained in the 
set $\Cal C$ of subcurves $Y\subseteq X(s)$ such that $\Cal I(s)$ 
decomposes at 
$Y$. If $\Cal C^1$ (or $\Cal C^{-1}$) is strictly contained in $\Cal C$, 
then we replace $\Cal I$ by $\Cal I^1$ (or $\Cal I^{-1}$) and start the above 
procedure again, but now with a ``better'' extension. If not, then both 
$\Cal I^1(s)$ and $\Cal I^{-1}(s)$ decompose at $Z$. In this case, let
$$
\Cal I^2:=\text{ker}(\Cal I^1 @>>> \Cal I^1(s)_Z) \text{\  \  and \  \  }
\Cal I^{-2}:=\text{ker}(\Cal I^{-1} @>>> \Cal I^{-1}_{Z^c}),
$$
and apply the argument used above for $\Cal I^1$ and $\Cal I^{-1}$ to 
both $\Cal I^2$ and $\Cal I^{-2}$. Applying the above procedure repeatedly, 
it is clear that 
we either obtain an extension $\Cal I$ of $I_{\eta}$ that is 
relatively simple over 
$S$, or we end up with two infinite filtrations of a certain 
extension $\Cal I$,
$$
\dots\subseteq \Cal I^i\subseteq\dots\subseteq \Cal I^1 \subseteq 
\Cal I^0:=\Cal I,
$$
$$
\dots\subseteq \Cal I^{-i} \subseteq \dots\subseteq \Cal I^{-1}
\subseteq\Cal I^0:=\Cal I,
$$
with quotients
$$
\frac{\Cal I^i}{\Cal I^{i+1}}=\Cal I^i(s)_Z \text{\  \  and \  \  } 
\frac{\Cal I^{-i}}{\Cal I^{-i-1}}=\Cal I^{-i}(s)_{Z^c}
$$
for $i\geq 0$, where $Z\subsetneqq X(s)$ is a non-empty, proper 
subcurve such that 
$\Cal I^i(s)$ decomposes at $Z$ for every integer $i$. We will show 
by contradiction that the latter 
situation is not possible. We may assume that $R$ is complete. (If not, 
just extend the 
sheaves $\Cal I^i$ over the completion of $R$.) By Lemma 17, 
there are $S$-flat quotients $\Cal F$ and $\Cal G$ of $\Cal
I$ such that $\Cal F(s)=\Cal I(s)_Z$ and $\Cal
G(s)=\Cal I(s)_{Z^c}$. Consider the induced 
homomorphism, $\phi\:\Cal I @>>> \Cal
F\oplus\Cal G$. By assumption, $\phi(s)$ is an
isomorphism. Since being an isomorphism is an open
property, then $\phi$ is an isomorphism. Thus, $I_{\eta}\cong\Cal I(\eta)$ 
is not simple, a contradiction. The proof of (2) is complete.

We prove (3) now. By (1), we may pick 
an extension $\Cal I$ of $I_{\eta}$. Consider the infinite filtration,
$$
\dots\subseteq\Cal I^i\subseteq\dots\subseteq \Cal
I^1\subseteq\Cal I^0:=\Cal I,
$$
with quotients
$$
\frac{\Cal I^i}{\Cal I^{i+1}}=\Cal I^i(s)_{Z_i},
$$
where $Z_i\subseteq X(s)$ 
is the maximum subcurve among the subcurves $W\subseteq 
X(s)$ with minimum $\beta_{\Cal I^i(s)}(W)$, for each $i\geq 0$. 
We claim that $\Cal I^i(s)$ is semistable with respect to $\Cal E(s)$ 
for some $i\geq 0$. 
Suppose by contradiction that our claim is false. We
may assume that $R$ is complete. (If not, just extend $\Cal E$ and the 
sheaves $\Cal I^i$ over the completion of $R$.) Since $\Cal I^i(s)$ is 
not semistable, then $Z_i$ 
is a non-empty, proper subcurve of $X(s)$ with
$\beta_{\Cal I^i(s)}(Z_i)<0$, for every $i\geq 0$. 
By Lemma 23, we may assume
that both $Z_i$ and $\beta_{\Cal I^i(s)}(Z_i)$ do not depend on $i$, 
and $\Cal I^i(s)$ decomposes at $Z_i$ for every $i\geq 0$. 
Let $Z:=Z_i$ and $\beta:=\beta_{\Cal I^i(s)}(Z_i)$ for every $i\geq 0$. 
By Lemma 17, there is an $S$-flat quotient $\Cal F$ of $\Cal I$ such
that $\Cal F(s)=\Cal I(s)_Z$. Since 
$\Cal F$ is $S$-flat and 
$\chi(\Cal F(s)\otimes\Cal E(s))=r\beta<0$, then also 
$\chi(\Cal F(\eta)\otimes\Cal E(\eta))<0$. Thus, 
$I_{\eta}\cong\Cal I(\eta)$ is 
not semistable with respect to $\Cal E(\eta)$. This contradiction 
shows that there is an extension 
$\Cal I$ of $I_{\eta}$ that 
is relatively semistable with respect to $\Cal E$. Suppose now that 
$I_{\eta}$ is simple. It is not necessarily true that 
$\Cal I$ is relatively 
simple. Nevertheless, we can apply the construction in the 
proof of (2) to 
$\Cal I$ to get a relatively simple sheaf that will still be 
relatively semistable with respect to $\Cal E$, by Lemma 18. 
The proof of (3) is complete.

We prove (4) now. By (3), there is a relatively 
semistable sheaf $\Cal I$ on $X$ over $S$ with respect to $\Cal E$ 
such that $\Cal I(\eta)\cong
I_{\eta}$. Consider the infinite filtration,
$$
\dots\subseteq\Cal I^i\subseteq\dots\subseteq\Cal
I^1\subseteq \Cal I^0:=\Cal I,
$$
with quotients
$$
\frac{\Cal I^i}{\Cal I^{i+1}}=\Cal I^i(s)_{Z_i},
$$
where $Z_i\subseteq X(s)$ is the minimum subcurve of 
$X(s)$ containing $\sigma(s)$ 
such that $\beta_{\Cal I^i(s)}(Z_i)=0$, for $i\geq 0$. We claim 
that $\Cal I^i(s)$ is $\sigma(s)$-quasistable 
with respect to $\Cal E(s)$ for 
some $i\geq 0$. In fact, it follows from Lemma 24 that
$$
Z_0\subseteq Z_1\subseteq\dots\subseteq Z_i\subseteq\dots.
$$
Thus, we may assume that $Z_i$ does not depend on $i$. Let 
$Z:=Z_i$ for every $i\geq 0$. We will show that $Z=X(s)$. 
It follows from Lemma 24 that $\Cal I^i(s)$ 
decomposes at $Z$ for every $i\geq 0$. 
We may now assume that $R$ is complete. (If not, 
just extend $\sigma$, $\Cal E$ and the sheaves $\Cal I^i$ over the 
completion of $R$.) By Lemma 17, there is an $S$-flat quotient $\Cal
F$ of $\Cal I$ such that $\Cal F(s)=\Cal
I(s)_Z$. Since $\beta_{\Cal I(s)}(Z)=0$, and
$\Cal F$ is $S$-flat, then $\chi(\Cal F(\eta)\otimes
\Cal E(\eta))=0$. Let $Y\subseteq X(\eta)$ 
be the maximum subcurve contained in the support of $\Cal
F(\eta)$. Since $\chi(\Cal F(\eta)\otimes\Cal E(\eta))=0$, and $\Cal I(\eta)$ 
is semistable, then $\Cal F(\eta)=\Cal I(\eta)_Y$ and 
$\beta_{\Cal I(\eta)}(Y)=0$. It follows that $\Cal F$ is relatively 
torsion-free on $X$ over $S$. So, $\sigma^*\Cal F$ is free. 
Since $\sigma^*\Cal F(s)\neq 0$, then also $\sigma^*\Cal F(\eta)\neq 0$. 
Thus, $Y$ contains $\sigma(\eta)$. 
Since $\Cal I(\eta)\cong I_{\eta}$ is $\sigma(\eta)$-quasistable, 
then $Y=X(\eta)$. It follows that $\Cal F(\eta)=\Cal I(\eta)$, and thus 
$\Cal F(s)=\Cal I(s)$. So, $Z=X(s)$, and thus $\Cal I(s)$ is 
$\sigma(s)$-quasistable. The proof of (4) is complete.

The proof of (5) will be
left to the reader. Roughly speaking, 
the proof consists of applying $n$ times 
the argument in the proof of (4), where 
$n$ is the number of irreducible components of $X(s)$.\qed
\enddemo

\vskip0.4cm

\parindent=0pt

\sl Remark 26. \rm Thm. 25 is not purely existencial. In fact, we 
have established a method to produce an extension of 
$I_{\eta}$ with the same ``good'' properties (semistability, 
quasistability, etc.) as $I_{\eta}$, given any 
extension $\Cal I$. We just construct a filtration,
$$
\dots\subseteq \Cal I^i \subseteq \dots\subseteq \Cal I^1 \subseteq 
\Cal I^0=\Cal I,
$$
of $\Cal I$ with quotients of the form:
$$
\frac{\Cal I^i}{\Cal I^{i+1}}=\Cal I^i(s)_{Z_i},
$$
where $Z_i\subseteq X(s)$ is a suitably chosen subcurve for each $i\geq 0$, 
as described in the proof of Thm. 25. Then, $\Cal I^i$ will 
be a ``good'' extension of $I_{\eta}$ for some $i\geq 0$. Note however 
that the minimum 
such $i$ depends on the original extension $\Cal I$.

\vskip0.8cm

\parindent=12pt

\heading
3. The fine moduli spaces
\endheading

Let $f\: X @>>> S$ be a flat, projective morphism whose
geometric fibres are curves. 
Let $\text{\bf J}^*$ denote 
the contravariant functor from the category of locally 
Noetherian $S$-schemes to sets, defined on an $S$-scheme $T$ by:
$$
\text{\bf J}^*(T):=\{\text{relatively simple,
torsion-free, rank }1 \text{ sheaves on } X\times_S T \text{ over } T\}/\sim,
$$
where ``$\sim$'' is the following equivalence relation:
$$
\Cal I_1\sim \Cal I_2 \Longleftrightarrow \text{ There is an
invertible sheaf } M \text{ on } T \text{ such that }
\Cal I_1\cong \Cal I_2\otimes M.
$$
Let $\text{\bf J}$ be the \'etale sheaf associated to
$\text{\bf J}^*$. By \cite{\bf 4\rm , Thm. 7.4, p. 99}, the
functor $\text{\bf J}$ is represented by an 
algebraic space $J$, locally of finite type 
over $S$. Note that the formation of $J$ commutes with
base change.

For every integer $d$, let $J_d\subseteq J$ be
the subspace parametrizing relatively simple, 
torsion-free, rank 1 sheaves with relative Euler
characteristic $d$. It is clear that $J_d$ is an open
subspace of $J$, and that $J$ is the disjoint
union of the $J_d$, for $d$ ranging
through all the integers. The formation of $J_d$ 
commutes also with base change.

Fix an integer $d$. Fix a vector
bundle $\Cal E$ on $X$ of rank $r>0$ and relative degree $-rd$ over $S$. 
We consider $\Cal E$ our relative polarization on $X$ over $S$. 
Let $J_{\Cal E}^s$ (resp. $J_{\Cal E}^{ss}$, resp. $J_{\Cal E}^{qs}$) 
denote the subspace 
of $J$ parametrizing relatively simple, torsion-free, rank 1 
sheaves on $X$ over $S$ that are relatively stable (resp.
semistable, resp. quasistable) with respect to $\Cal E$. 
If $\sigma\: S @>>> X$ is a section of $f$ through the $S$-smooth 
locus of $X$, we let $J_{\Cal E}^{\sigma}$ denote the subspace of $J$ 
parametrizing relatively simple, 
torsion-free, rank 1 sheaves on $X$ over $S$ that are 
relatively $\sigma$-quasistable with respect to $\Cal E$ over $S$. 
It is clear from the definitions in Sect. 1 that
$$
J_{\Cal E}^s\subseteq J_{\Cal E}^{\sigma}\subseteq 
J_{\Cal E}^{qs}\subseteq J_{\Cal E}^{ss}\subseteq J_d.
$$
The formations of all the above spaces commute with base change.

\vskip0.4cm

\proclaim{Proposition 27} The subspaces $J_{\Cal E}^s,J_{\Cal E}^{\sigma},
J_{\Cal E}^{qs},J_{\Cal E}^{ss}\subseteq J_d$ are open.
\endproclaim

\demo{Proof} The proposition follows by semicontinuity 
from the cohomological characterizations 
of Subsect. 1.4, namely Thm. 8, Thm. 9 and Cor. 10, modified 
according to Rmk. 15. We leave it to the reader to work out the details.\qed
\enddemo

\vskip0.4cm

\parindent=0pt

\sl Proof of Theorem A. \rm We may assume that $S$ is Noetherian. It follows 
straight from the definitions that the family of all relatively 
semistable sheaves on $X$ over $S$ is bounded. Hence, 
$J^{ss}_{\Cal E}$ is of finite type over $S$. 
Since $J^{ss}_{\Cal E}$ is of finite type over $S$, then 
statement (1) follows from
Thm. 25, statement (2) follows from Prop. 19, and statement (3) 
follows from Prop. 20 and Thm. 25. The proof is complete.\qed

\vskip0.4cm

\parindent=12pt

If $Y$ is a curve over a field $k$, we let $\overline{Y}:=Y\times_k \bar k$, 
where $\bar k$ denotes the algebraic closure of $k$.

\vskip0.4cm

\proclaim{Lemma 28} Let $s_0\in S$ be a closed point. 
Let $n$ denote the number of 
irreducible components of the geometric fibre $\overline{X(s_0)}$. 
Then, there are an \'etale morphism, 
$S' @>>> S$, containing $s_0$ in its image, and sections
$$
\sigma_1,\dots,\sigma_n\: S' @>>> X\times_S S'
$$
of the projection morphism $X\times_S S' @>>> S'$, such that:
\roster
\item $\sigma_i$ factors through the $S'$-smooth locus of $X\times_S S'$ for 
$i=1,\dots,n$;
\item for every $s\in S'$, every irreducible component of $X(s)$ is 
geometrically integral, and contains $\sigma_i(s)$ for some $i$.
\endroster
\endproclaim

\demo{Proof} Replacing $S$ by an \'etale
neighbourhood of $s_0$, if necessary, we may assume that the
irreducible components, $X_1,\dots,X_n$, of $X(s_0)$ are geometrically
integral. For each $i=1,\dots,n$, let
$$
Y_i:=X\setminus (\bigcup_{j\neq i} X_j).
$$
By \cite{\bf 16\rm , IV-4-17.16.3}, for each $i=1,\dots, n$ there are 
an \'etale morphism, $S_i @>>> S$, containing $s_0$ in its image, 
and an $S$-morphism, $\sigma_i\: S_i @>>> X$, 
factoring through the $S$-smooth locus of $Y_i$. We may assume that 
$S=S_1=\dots=S_n$. By construction, $\sigma_i(s_0)\in X_i$ for 
$i=1,\dots,n$. Thus, $\Cal
O_{X(s_0)}(\sigma_1(s_0)+\dots+\sigma_n(s_0))$ is ample on $X(s_0)$. 
Since ampleness is an open property \cite{\bf 16\rm , 
III-1-4.7.1, p. 145}, replacing $S$ by a Zariski neighbourhood of $s$, 
if necessary, we may assume
that the sheaf $\Cal O_{X(s)}(\sigma_1(s)+\dots+\sigma_n(s))$ is
ample for every $s\in S$. Consequently, for every $s\in S$, every 
irreducible component of $X(s)$ 
is geometrically integral, and contains 
$\sigma_i(s)$ for some $i$. The proof is complete.\qed
\enddemo

\vskip0.4cm

\proclaim{Lemma 29} Assume that there are sections 
$\sigma_1,\dots,\sigma_n\: S @>>> X$ of $f$ such that:
\roster
\item $\sigma_i$ factors through the $S$-smooth locus of $X$ for 
$i=1,\dots,n$;
\item for every $s\in S$, every irreducible component of $X(s)$ is 
geometrically integral, and contains $\sigma_i(s)$ for some $i$.
\endroster 
Then, there are relative polarizations $\Cal E_j$ on $X$ over $S$ such that
$$
J=\bigcup_j J^s_{\Cal E_j}.
$$
\endproclaim

\demo{Proof} Let $\underline e:=(e_1,\dots,e_n)$ be a $n$-uple of integers. 
For $i=1,\dots,n$, let
$$
\Cal L_i:=\Cal O_X(-e_i\Sigma_i),
$$
where $\Sigma_i$ is the relative effective Cartier divisor on $X$ 
over $S$ corresponding to $\sigma_i$. 
Let $r$ be a positive integer. Let 
$$
\Cal E_{(\underline e,r)}:=\Cal L_1\oplus\dots\oplus\Cal L_n
\oplus\Cal O_X^{\oplus n(r-1)}.
$$
It follows from Cor. 13 that
$$
J=\bigcup_{(\underline e,r)} J^s_{\Cal E_{(\underline e,r)}},
$$
where the union runs through all integers 
$r\geq 1$ and all $n$-uples $\underline e$ such that $nr|(e_1+\dots+e_n)$. 
The proof is complete.\qed
\enddemo

\vskip0.4cm

\parindent=0pt

\sl Remark 30. \rm Let $\Cal M$ be an invertible sheaf on $X$. 
Let $s\in S$, and $F$ be a vector 
bundle on $X(s)$ with $\det F\cong\Cal M(s)$. 
Then, there are an \'etale (Zariski, if $k(s)$ is infinite) covering 
$S' @>>> S$, a point $s'\in S'$ lying above $s$, and a 
vector bundle $\Cal F'$ on $X\times_S S'$ such that:
\roster
\item $\Cal F'(s')\cong F\otimes_{k(s)} k(s')$;
\item $\det\Cal F'\cong\Cal M\otimes\Cal O_{S'}$.
\endroster 
The proof of the above fact can be 
carried out using an argument similar to that in Rmk. 15.

\parindent=12pt

\vskip0.4cm

\proclaim{Lemma 31} Let $\Sigma\subseteq J$ be an open subspace of 
finite type over $S$. Then, there is an \'etale covering
$S'$ of $S$ such that $\Sigma\times_S S'$ is a scheme.
\endproclaim

\demo{Proof} The statement is clearly local on $S$ for
the \'etale topology. We may thus assume that $S$ is Noetherian. 
Moreover, in view of Lemma 28 and Lemma 29, we may assume that there 
are a section, $\sigma\: S @>>> X$, of $f$ through the $S$-smooth locus of 
$X$, and a relative polarization $\Cal E$ on $X$ over $S$, and we need only 
show that there is an \'etale covering $S'$ of $S$ such that 
$J^{\sigma}_{\Cal E}\times_S S'$ is a scheme.

Fix $s_0\in S$. We will show that there is an \'etale
neighbourhood $S'$ of $s_0$ in $S$ such that $J^{\sigma}_{\Cal E}\times_S
S'$ is a scheme. It follows from Thm. 9 (adequately 
modified by Rmk. 15), Rmk. 30, and the fact 
that $J^{\sigma}_{\Cal E}$ is Noetherian, that, up to replacing $S$ by an 
\'etale neighbourhood of $s_0$ in $S$, if necessary, 
there are finitely many vector 
bundles $\Cal F_1,\dots,\Cal F_t$ on $X$ such that
$$
J^{\sigma}_{\Cal E}(s_0)=\bigcup_{i=1}^t U_i(s_0),
$$
where $U_i\subseteq J_{\Cal E}^{\sigma}$ 
is the open
subspace parametrizing the relatively $\sigma$-quasistable sheaves 
$\Cal I$ on $X$ over $S$ 
with respect to $\Cal E$ such that 
$$
h^0(X(s),\Cal I(s)\otimes\Cal F_i(s))=0 \text{\  \  and \  \  } 
h^1(X(s),\Cal I(s)\otimes\Cal F_i(s))=1,
$$
and the unique (modulo $k(s)^*$) non-zero homomorphism 
$\Cal I(s) @>>> \Cal F_i(s)^*\otimes\omega(s)$ 
is injective for every $s\in S$. As in 
\cite{\bf 12\rm , Cor. 2}, it is clear that the space $U_i$ is 
isomorphic to a certain open subscheme 
of $\text{Quot}_{X/S}(\Cal F_i^*\otimes\omega)$, Grothendieck's 
scheme of quotients of $\Cal F_i^*\otimes\omega$. In particular, 
$U_i$ is a scheme. 
Since $J_{\Cal E}^{\sigma}$ is universally closed over $S$ by 
Thm. A, then there is an open neighbourhood $S'\subseteq S$ of 
$s_0$ such that
$$
J_{\Cal E}^{\sigma}\times_S S'=\bigcup_{i=1}^t U_i\times_S S'.
$$
Thus, $J_{\Cal E}^{\sigma}\times_S S'$ is a scheme. 
The proof is complete.\qed
\enddemo

\vskip0.4cm

In Sect. 5 we will drop the hypothesis that $\Sigma$ is of finite type 
over $S$. More precisely, Cor. 42 asserts the conclusion of 
Lemma 31 for $J$, instead of just $\Sigma$.

\vskip0.8cm

\heading
4. Certain special cases
\endheading

Let $X$ be a connected, reduced curve over an 
algebraically closed field $k$. Let $g$ denote the arithmetic genus of 
$X$. Let $J$ denote the algebraic space parametrizing simple, 
torsion-free, rank 1 sheaves on $X$. 
It follows from Lemma 29 and Lemma 31 that $J$ is a scheme, locally 
of finite type over $k$.

\vskip0.4cm

\parindent=0pt

\sl Example 32. \rm (Joining two curves) 
Assume that there are subcurves $Y,Z\subseteq X$ covering $X$ such 
that $Y\cap Z$ has length 1. Then, 
every simple, torsion-free, rank 1 
sheaf on $X$ must be invertible on $Y\cap Z$. Let $J_X$ (resp. $J_Y$, 
resp. $J_Z$) denote the algebraic space parametrizing simple, torsion-free, 
rank 1 sheaves on $X$ (resp. $Y$, resp. $Z$). It follows from 
our above observation that we have a morphism,
$$
J_X @>>> J_Y\times J_Z,
$$
defined by restriction of sheaves on $X$ to $Y$ and $Z$. We can 
show that the above morphism is an isomorphism. 

\vskip0.4cm

\parindent=0pt

\sl Example 33. \rm (Abel maps) Let 
$$
\delta_X:=\min_{Y\subsetneqq X} \chi(\Cal O_{Y\cap Y^c}),
$$
where $Y$ runs through all non-empty, proper subcurves of $X$. 
(If $X$ is irreducible, let $\delta_X:=\infty$.) Since 
$X$ is connected, then $\delta_X>0$. If $\delta_X=1$, then we are in the 
case of Ex. 32. Let 
$I\subseteq\Cal O_X$ be the ideal sheaf of a subscheme $D\subset X$ 
of finite length $m$. 
Let $Y\subsetneqq X$ be a non-empty, proper subcurve, and let $Z:=Y^c$. 
Then, there is a natural commutative diagram,
$$
\CD 
0 @>>> I @>>> \Cal O_X @>>> \Cal O_D @>>> 0\\
@. @VVV @VVV @VVV @.\\
0 @>>> I_Y\oplus I_Z @>>> \Cal O_Y\oplus\Cal O_Z @>>> \Cal O_{D\cap Y}
\oplus\Cal O_{D\cap Z} @>>> 0,
\endCD
$$
with exact horizontal sequences. It follows that
$$
\chi(I_Y)+\chi(I_Z)-\chi(I)\geq\chi(\Cal O_{Y\cap Z})-m,
$$
with equality if and only if $D\subseteq Y\cap Z$. Therefore, if 
$m<\delta_X$, then $I$ is simple. On the other hand, if 
$D=Y\cap Y^c$, for any non-empty, proper subcurve 
$Y\subsetneqq X$, then $I$ is not simple. To summarize, 
there are ideal sheaves of 
subschemes $D\subset X$ of length $\delta_X$ which are not simple, and 
$\delta_X$ is the minimum length where such phenomenon occurs.

\parindent=12pt

Let $m$ be any 
non-negative integer with $m<\delta_X$. 
Let $H_m$ denote the Hilbert scheme of $X$, 
parametrizing length $m$ subschemes of $X$. 
Of course, $H_1=X$. Let $M$ be an invertible sheaf on $X$. 
For every subscheme $D\subseteq X$, we let 
$\Cal M_D\subseteq\Cal O_X$ denote its ideal sheaf. 
It follows from our above considerations that we have a well-defined 
morphism,
$$
\aligned 
\Gamma^m_M\: H_m @>>> & \  \  \  \  \  J_d\\
[D]\  \mapsto & [\Cal M_D\otimes M],
\endaligned
$$ 
where $d:=\chi(M)-m$. 
The morphism $\Gamma^m_M$ is called the \sl Abel map in degree $m$ \rm of 
$X$. If $\delta_X>1$, we let
$$
\Gamma_M:=\Gamma_M^1\: X @>>> J_d.
$$
If $X$ is irreducible, then $\Gamma_M$ is a closed embedding 
\cite{\bf 4\rm , Thm. 8.8, p. 108}. If, in addition, $g=1$, then 
$\Gamma_M$ is an isomorphism \cite{\bf 4\rm , Ex. 8.9, p. 109}. In 
general, I do not know what $\Gamma_M$ is, but we will treat the case 
$g=1$ below.

\vskip0.4cm

\parindent=0pt

\sl Example 34. \rm (Genus 1 curves) Assume that $g=1$ and $\delta_X>1$. 
We claim that $\Cal O_X$ is the dualizing sheaf of 
$X$. In fact, there is a non-zero section, $h\:\Cal O_X @>>> \omega$, 
where $\omega$ is the dualizing sheaf of $X$. If we proved that 
$h$ is injective, then it would follow from degree considerations 
that $h$ is an isomorphism. 
Let $Z\subseteq X$ be 
the non-empty subcurve such that $\Cal O_Z\cong\text{im}(h)$. Then, 
$h$ factors through $\omega_Z\subseteq\omega$, where $\omega_Z$ is the 
dualizing sheaf of $Z$. Since $X$ has genus 1, and 
$\Cal O_Z\hookrightarrow\omega_Z$ is injective, 
then it follows that $Z$ is connected and 
has arithmetic genus 1. To prove our claim, we need only show that 
$X$ does not admit any proper, connected subcurve $Z\subsetneqq X$ 
of arithmetic genus 1. Suppose, by contradiction, that $Z$ were such a curve. 
Then, $Y:=Z^c$ has arithmetic genus 0, 
and $h^0(Y,\Cal O_Y)=\chi(\Cal O_{Y\cap Z})$. 
It follows that $Y$ has exactly $\chi(\Cal O_{Y\cap Z})$ connected 
components. Then, each connected component $Y_0\subseteq Y$ 
intersects $Z$ in a scheme of length 1. Thus, 
$\chi(\Cal O_{Y_0\cap Y_0^c})=1$, contradicting the fact that 
$\delta_X>1$. 

\parindent=12pt

Let $d$ be an integer. Let $M$ be an invertible sheaf on $X$ 
of degree $d+1$. Let $J_M\subseteq J_d$ denote the subset 
parametrizing simple, torsion-free, rank 1 sheaves $I$ 
on $X$ such that $\chi(I_Y)\geq \deg_Y M$ for every non-empty, 
proper subcurve $Y\subsetneqq X$. We claim that $J_M\subseteq J_d$ is a 
complete, open subscheme, and $\Gamma_M$ factors through $J_M$.

First, let $p\in X$ be any non-singular point, and put 
$E:=M^*\otimes\Cal O_X(p)$. It follows easily from the definition that 
a torsion-free, rank 1 sheaf $I$ on $X$ with $\chi(I)=d$ is $p$-quasistable 
with respect to $E$ if and only if $\chi(I_Y)\geq \deg_Y M$ 
for every non-empty, proper subcurve $Y\subsetneqq X$. Thus, 
$J_M=J^p_E$. It follows from Thm. A and Prop. 27 that 
$J_M$ is a complete, open subscheme of $J_d$.

Second, let $q\in X$ and $Y\subsetneqq X$ be a connected, proper subcurve. 
Since $Y$ has arithmetic genus 0, then $\chi(\Cal M_{q,Y})=1$ if 
$q\not\in Y$, and $\chi(\Cal M_{q,Y})=0$ otherwise. At any rate, we 
obtain that $\chi(I_Y)\geq \deg_Y M$, where $I:=\Cal M_q\otimes M$. The 
upshot is that $\Gamma_M$ factors through $J_M$, proving our claim.

We claim now that $\Gamma_M$ is an isomorphism onto $J_M$. 
To show our claim, we will construct the inverse morphism, 
$\Lambda_M\: J_M @>>> X$, as follows: 
Let $I$ be a simple, torsion-free, rank 1 sheaf on $X$ such that $\chi(I)=d$, 
and
$$
\chi(I_Y)\geq \deg_Y M\tag{34.1}
$$
for every non-empty, proper subcurve 
$Y\subsetneqq X$. We will show that 
$$
H^0(X,I\otimes M^*)=0,\tag{34.2}
$$
and the unique (modulo $k^*$) non-zero homomorphism, 
$\lambda\: I\otimes M^* @>>> \Cal O_X$, is an isomorphism onto the ideal 
sheaf $\Cal M_q$ of a point $q\in X$. We will then let 
$\Lambda_M([I]):=q$. It is clear that $\Lambda_M$, if defined, is 
inverse to $\Gamma_M$.

In order to show (34.2), 
suppose by contradiction that there is a non-zero section, $\mu\: 
\Cal O_X @>>> I\otimes M^*$. Let $Z\subseteq X$ be the subcurve such that 
$\Cal O_Z\cong\text{im}(\mu)$. Then, $\mu$ factors through 
$J\otimes M^*$, where $J:=\text{ker}(I \twoheadrightarrow I_{Z^c})$. 
It follows from (34.1), and the fact that $\chi(I\otimes M^*)=-1$, 
that $\chi(J\otimes M^*)\leq -1$. 
Since $\Cal O_Z \hookrightarrow J\otimes M^*$, then it follows that 
$\chi(\Cal O_Z)\leq -1$. But, since 
$Z\subseteq X$, and $g=1$, then $\chi(\Cal O_Z)\geq 0$. This 
contradiction proves (34.2).

Since $\Cal O_X$ is the dualizing sheaf of $X$, then it follows from (34.2), 
and the fact that $\chi(I\otimes M^*)=-1$, that 
there is a unique (modulo $k^*$) non-zero homomorphism 
$\lambda\: I\otimes M^* @>>> \Cal O_X$. We need only show that 
$\lambda$ is injective. Suppose by contradiction 
that $\lambda$ is not injective.
Let $Z\subseteq X$ be the subcurve 
of $X$ such that $I_Z\otimes M^*\cong\text{im}(\lambda)$. By contradiction 
hypothesis, we have that $Z\subsetneqq X$ is a non-empty, proper subcurve. 
It follows from (34.1) that $\chi(I_Z\otimes M^*)\geq 0$. Let $J$ denote the 
ideal sheaf of $Z^c$. Since $Z^c$ has arithmetic genus 0, then 
$\chi(J)=-1$. But, since $I_Z\otimes M^*\hookrightarrow J$ is injective, 
then $\chi(J)\geq\chi(I_Z\otimes M^*)$. We have a contradiction, proving that 
$\lambda$ is injective. 

We observe that we have naively defined $\Lambda_M$ as a set map, 
but it is clearly possible to 
apply the above argument to a family of torsion-free, rank 1 
sheaves on $X$, and thus define $\Lambda_M$ as a morphism.

\vskip0.4cm

\parindent=0pt

\sl Example 35. \rm (Locally planar curves) Assume that $X$ is a 
locally planar curve. In other words, assume that $X$ can be embedded into 
a smooth surface. Then, Altman, Kleiman and Iarrobino 
\cite{\bf 2\rm , Cor. 7, p. 7} showed that the 
Hilbert scheme $H_m$ of $X$, parametrizing length $m$ subschemes of $X$, 
is $m$-dimensional, reduced and a local complete intersection. In addition, 
the proof 
of \cite{\bf 2\rm , Thm. 9, p. 8} can be easily adapted to show 
that $J$ is reduced, a local complete intersection, and has pure 
dimension $g$ at every point. 
Finally, it follows from \cite{\bf 2\rm , Thm. 5, p. 5} that 
the invertible sheaves form an open dense subscheme of $J$.

\vskip0.4cm

\parindent=0pt

\sl Example 36 \rm (Two-component curves) Assume that 
$X$ has only two irreducible components, $X_1$ and $X_2$. Let 
$\delta$ denote the length of $X_1\cap X_2$. Let $E$ be a polarization on 
$X$ of rank $r$ and degree $-rd$, for an integer $d$. 
Let $e_i:=-\deg_{X_i} E$ for $i=1,2$. 
For $i=1,2$, let $J_E^i$ denote the moduli space of $X_i$-quasistable sheaves 
on $X$ with respect to $E$. We have two cases:
\roster 
\item $r\not|e_1$: In this case, 
$J_E^s=J_E^{ss}$, and $J_E^s$ is complete. If 
$X$ is locally planar, then $J_E^s$ has $\delta$ irreducible 
components.
\item $r|e_1$: In this case, 
$$
J_E^s\subsetneqq J_E^1, J_E^2\subsetneqq J_E^{ss}.
$$
If $X$ is locally planar, then 
$J_E^{ss}$ has $\delta+1$ irreducible components, whereas 
$J_E^1$ and $J_E^2$ have $\delta$ components each, 
and $J_E^s$ has $\delta-1$ components. 
As we had already observed in Ex. 5, we have that 
$J_E^{qs}=J_E^{ss}$.
\endroster
Case 1 corresponds to Caporaso's general case \cite{\bf 10\rm , 7.3, p. 646}, 
whereas Case 2 corresponds to her special case. 

\parindent=12pt
 
Assume now that $X_1$ and $X_2$ are smooth, and 
intersect at two ordinary nodes. 
So, $\delta_X=2$, and $X$ is locally planar. Denote by $P_{d+1}$ the 
Jacobian of $X$, parametrizing invertible sheaves of Euler characteristic 
$d+1$, and consider the (well-defined) morphism:
$$
\aligned 
\Gamma\: X\times P_{d+1} @>>> & \  \  \  \  \  J_d\\
(q,M)\  \mapsto & [\Cal M_q\otimes M].
\endaligned
$$
It is not difficult to show that $\Gamma$ is surjective, and 
smooth with relative dimension 1. 
In Case 1, we have that $J_E^s$ is the image under $\Gamma$ 
of a connected component of $X\times P_{d+1}$. 
In Case 2, we have that $J_E^{ss}=J_E^1\cup J_E^2$, and both $J_E^1$ and 
$J_E^2$ are images under $\Gamma$ of different connected components of 
$X\times P_{d+1}$. In fact, we have that $J_E^i=\Gamma(X\times P_E^i)$, 
where $P_E^i$ is the connected component of 
$P_{d+1}$ parametrizing 
invertible sheaves on $X$ with Euler characteristic $e_i/r+2$ 
on $X_i$, for $i=1,2$. The patching of $J_E^1$ and $J_E^2$ to produce 
$J_E^{ss}$ occurs on $J_E^s$, which is the image under $\Gamma$ 
of both $(X_1\setminus X_2)\times P_E^1$ and $(X_2\setminus X_1)\times P_E^2$. 

We observe the difference between our description and 
Caporaso's in loc. cit.. 
As we are dealing with fine moduli spaces, without 
making identifications to guarantee separatedness, our moduli spaces 
$J_E^{ss}$ in Case 2 have more components than Caporaso's. In fact, to reach 
Caporaso's spaces, we have to identify and collapse 
the two non-stable components of $J_E^{ss}$ 
through the JH-equivalence into a positive codimensional locus. Such 
identification will be carried out in the next section, by means of theta 
functions.

\vskip0.8cm

\parindent=12pt

\heading
5. Theta functions and the separated realizations
\endheading

\parindent=0pt

\subheading{5.1. The determinant of cohomology} 
Let $f\: X @>>> S$ be a flat, projective morphism whose
geometric fibres are curves. 
Let $\Cal F$ be an $S$-flat coherent sheaf on $X$. 
The \sl determinant of cohomology of $\Cal F$ with respect
to $f$ \rm is the invertible sheaf $\Cal D_f(\Cal F)$ on
$S$ constructed as follows: Locally on $S$ there is a
complex
$$
0 @>>> G^0 @>\lambda>> G^1 @>>> 0
$$
of free sheaves of finite rank such that, for every coherent sheaf $M$
on $S$, the cohomology groups of $G^{\bullet}\otimes M$
are equal to the higher direct images of $\Cal F\otimes M$
under $f$. The complex $G^{\bullet}$ is unique up to
unique quasi-isomorphism. Hence, its determinant,
$$
\det G^{\bullet}:=(\bigwedge^{\text{rank } G^1}
G^1)\otimes(\bigwedge^{\text{rank } G^0} G^0)^{-1},
$$
is unique up to canonical isomorphism. 
The uniqueness allows us to glue together the
local determinants to obtain the invertible sheaf $\Cal
D_f(\Cal F)$ on $S$. 

\parindent=12pt

The determinant of cohomology has the following properties:
\roster
\item \sl Functorial property\rm : We can functorially associate to any 
isomorphism, $\phi\: \Cal F_1 \cong \Cal F_2$, of $S$-flat coherent sheaves 
on $X$ an isomorphism: 
$$
\Cal D_f(\phi)\: \Cal D_f(\Cal F_1) \cong \Cal D_f(\Cal F_2).
$$ 
\item \sl Additive property\rm : We can functorially associate to 
any short exact sequence,
$$
\alpha\:\  0 @>>> \Cal F_1 @>>> \Cal F_2 @>>> \Cal F_3 @>>> 0,
$$
of $S$-flat coherent sheaves on $X$ an isomorphism:
$$
\Cal D_f(\alpha)\: \Cal D_f(\Cal F_2) \cong 
\Cal D_f(\Cal F_1)\otimes\Cal D_f(\Cal
F_3).
$$
\item \sl Projection property\rm : We can functorially associate to 
any $S$-flat coherent sheaf $\Cal F$ on $X$ of relative 
Euler characteristic $d$ over $S$, and 
any invertible sheaf $M$ on $S$, an isomorphism:
$$
\Cal D_f(\Cal F\otimes M)\cong\Cal D_f(\Cal
F)\otimes M^{\otimes d}.
$$
\item \sl Base-change property\rm : We can functorially associate 
to any $S$-flat coherent sheaf $\Cal F$ on $X$, and any Cartesian 
diagram of the form:
$$
\CD
X_1 @>h_1>> X\\
@Vf_1VV @VfVV\\
S_1 @>h>> S,
\endCD
$$
a base-change isomorphism:
$$
h^*\Cal D_f(\Cal F)\cong \Cal D_{f_1}(h_1^*\Cal F).
$$
\endroster 

For a more systematic development of the theory of
determinants, see \cite{\bf 20}. It is also possible to adopt
a more concrete approach to define $\Cal D_f(\Cal F)$, like the one 
used in \cite{\bf 6\rm , Ch. IV, \S 3}. 

If $\Cal F$ is an $S$-flat coherent sheaf on $X$ 
of relative Euler characteristic 0 over $S$, then
there is a canonical global section $\sigma_{\Cal F}$ of $\Cal
D_f(\Cal F)$, constructed as follows: Since 
$\chi(\Cal F(s))=0$ for every $s\in S$, then the ranks of
$G^0$ and $G^1$ in the local complex $G^{\bullet}$ are
equal. By taking the determinant of $\lambda$ we obtain a 
section of $\det G^{\bullet}$. Since the complex $G^{\bullet}$ is unique, 
such section is also unique, allowing us to glue the
local sections to obtain $\sigma_{\Cal F}$. 
The zero locus of $\sigma_{\Cal F}$ 
parametrizes the points $s\in S$ such that
$$
h^0(X(s),\Cal F(s))=h^1(X(s),\Cal F(s))\neq 0.
$$
(Another way of viewing $\sigma_{\Cal F}$ is as a generator of the 
0-th Fitting ideal of $R^1f_*\Cal F$.)

The global section $\sigma_{\Cal F}$ satisfies properties 
compatible with those of $\Cal D_f(\Cal F)$. For instance, 
if $\phi\: \Cal F_1 
\cong \Cal F_2$ is an isomorphism of $S$-flat coherent sheaves on $X$ with 
relative Euler characteristic 0 over $S$, then $\sigma_{\Cal F_2}=
\Cal D_f(\phi)\circ\sigma_{\Cal F_1}$. In addition, we have the 
additive property: If 
$$
\alpha\: 0 @>>> \Cal F_1 @>>> \Cal F_2 @>>> \Cal F_3 @>>> 0
$$
is a short exact sequence of $S$-flat coherent sheaves on $X$ with 
relative Euler characteristic $0$ over $S$, then
$$
\sigma_{\Cal F_1}\otimes\sigma_{\Cal F_3}=
\Cal D_f(\alpha)\circ\sigma_{\Cal F_2}.
$$
We leave it to the 
reader to state the projection and base-change properties 
of $\sigma_{\Cal F}$.

\vskip0.4cm

\parindent=0pt

\subheading{5.2. Theta functions} Fix an integer $d$, and let 
$\Sigma\subseteq J_d$ be an open subspace of finite type over $S$. 
Let $\Cal E$ be a vector bundle on $X$ of rank $r>0$ and relative degree 
$-rd$ over $S$. 

\parindent=12pt

(Notation: If $S^{(n)} @>>> S$ is a morphism of schemes, 
we let $X^{(n)}:=X\times_S S^{(n)}$ and 
$\Sigma^{(n)}:=\Sigma\times_S S^{(n)}$. We let $\Cal D^{(n)}$ denote the 
determinant of cohomology with respect to the projection morphism, 
$X^{(n)}\times_{S^{(n)}}\Sigma^{(n)} @>>> \Sigma^{(n)}$.)

Let $S' @>>> S$ be an \'etale covering such that 
$X'\times_{S'} \Sigma'$ admits a universal relatively torsion-free, 
rank 1 sheaf $\Cal I'$ over $\Sigma'$, and $\Sigma'$ is 
a scheme (see Lemma 31). Let $\Cal L'_{\Cal E}:=
\Cal D'(\Cal I'\otimes\Cal E)$. 
We can show, by faithfully flat descent, that the invertible sheaf 
$\Cal L'_{\Cal E}$ on $\Sigma'$ descends to a unique (up to isomorphism) sheaf 
$\Cal L_{\Cal E}$ on $\Sigma$. Since the 
procedure used is quite standard, we just list below the 
main properties that allow us to carry out descent, and mention 
how they are used:
\roster
\item Since $\Cal I'$ is universal and relatively simple over 
$\Sigma'$, then, for any $S'$-scheme $S''$ and any universal sheaf 
$\Cal I''$ on $X''\times_{S''}\Sigma''$, there are an invertible sheaf
$N$ on $\Sigma''$ and an isomorphism $\Cal I''\cong\Cal I'\otimes N$ 
on $X''\times_{S''} \Sigma''$ (see \cite{\bf 4\rm , Lemma 4.3, p. 79}).
\item Since the relative Euler characteristic of $\Cal I'\otimes
\Cal E$ over $\Sigma'$ is 0, then, by the projection and 
base-change properties of the determinant of cohomology, there are 
natural isomorphisms,
$$
\Cal D''(\Cal I'\otimes \Cal E\otimes N)\cong
\Cal D''(\Cal I'\otimes\Cal E\otimes\Cal O_{S''})
\cong \Cal D'(\Cal I'\otimes\Cal E)\otimes\Cal O_{S''},
$$
for any $S'$-scheme $S''$ and any 
invertible sheaf $N$ on 
$\Sigma''$.
\item Since $\Cal I'$ is relatively simple over 
$\Sigma'$, then, for any $S'$-scheme $S''$, and any invertible 
sheaf $N$ on $\Sigma''$, any 
isomorphism $\Cal I'\cong\Cal I'\otimes N$ 
is just an isomorphism $\Cal O_{\Sigma''}\cong N$ 
tensored with $\Cal I'$ \cite{\bf 5\rm , Lemma 5, p. 119}.
\endroster 
It follows from properties (1) and (2) 
that there are 
patching isomorphisms for $\Cal L'_{\Cal E}$ 
over the ``double intersection'' covering $S'':=S'\times_S S'$ of $S$. 
Property (3) 
implies 
that these patching homomorphisms satisfy the ``cocycle condition'' over the 
``triple intersection'' covering $S''':=S'\times_S S'\times_S S'$ of $S$. 
Thus, $\Cal L'_{\Cal E}$ descends to a sheaf $\Cal L_{\Cal E}$ on 
$\Sigma$. Properties (1), (2) and (3) show also that the descended sheaf, 
$\Cal L_{\Cal E}$, does not depend on the 
choice of an \'etale covering $S' @>>> S$, or of 
a particular universal sheaf $\Cal I'$ on 
$X'\times_{S'} \Sigma'$. The sheaf $\Cal L_{\Cal E}$ will be 
called the \sl Theta line bundle \rm on 
$\Sigma$ with respect to $\Cal E$. 

There is a natural way to 
produce sections of tensor powers of $\Cal L_{\Cal E}$, 
as we will show next.

Let $\Cal F$ be a vector bundle on $X$, with rank $mr$ and 
$\det \Cal F\cong(\det\Cal E)^{\otimes m}$ for some $m>0$. Assume 
that the determinant of cohomology of $\Cal F$ with respect to $f$ 
is equal to that of $\Cal E^{\oplus m}$. In other words, assume that 
there is an isomorphism 
$\Cal D_f(\Cal F)\cong\Cal D_f(\Cal E)^{\otimes m}$. 
Of course, such an isomorphism is defined modulo $H^0(S,\Cal O_S)^*$. 
Assume that $\Sigma$ is a scheme, and that there 
is a universal relatively torsion-free, rank 1 
sheaf $\Cal I$ on $X\times_S \Sigma$ 
over $\Sigma$. It 
follows from \cite{\bf 14\rm , Thm. I.1, p. 509} 
that there is a canonical, functorial isomorphism,
$$
\Cal D(\Cal I\otimes\Cal F)\otimes
\Cal D_f(\Cal E^{\oplus m})\cong
\Cal D(\Cal I\otimes\Cal E^{\oplus m})\otimes
\Cal D_f(\Cal F),
$$
that is well defined modulo $H^0(S,\Cal O_S)^*$. Using that 
$\Cal D_f(\Cal F)\cong\Cal D_f(\Cal E)^{\otimes m}$, we get a canonical 
(modulo $H^0(S,\Cal O_S)^*$) functorial isomorphism,
$$
\Phi_{\Cal F/\Cal E}\: \Cal D(\Cal I\otimes\Cal F)\cong
\Cal L_{\Cal E}^{\otimes m}.
$$
By means of $\Phi_{\Cal F/\Cal E}$, we may regard the global 
section $\sigma_{\Cal I\otimes\Cal F}$ of $\Cal D(\Cal I\otimes\Cal F)$ 
as a global section of $\Cal L_{\Cal E}^{\otimes m}$. 
More precisely, we let
$$
\theta_{\Cal F}:=\Phi_{\Cal F/\Cal E}\circ\sigma_{\Cal I\otimes\Cal F}\in 
H^0(\Sigma,\Cal L_{\Cal E}^{\otimes m}).
$$
The section $\theta_{\Cal F}$ is called the \sl theta function 
associated with $\Cal F$\rm . It is well defined modulo 
$H^0(S,\Cal O_S)^*$.

We observe that we don't really have to assume that $\Sigma$ is a scheme, or 
that $X\times_S \Sigma$ admits a universal sheaf over $\Sigma$, 
as we could have argued as before to obtain the section $\theta_{\Cal F}$ 
of $\Cal L_{\Cal E}^{\otimes m}$ by faithfully flat descent.

The theta functions will allow us to produce a separated 
realization of $\Sigma$, as we will show next.

\vskip0.4cm

\parindent=0pt

\subheading{5.3. The separated realizations} 
Let $\pi\: \Sigma @>>> S$ denote the structure morphism. Since 
$S$ is locally Noetherian, and $\pi$ is of finite type, then 
$\pi_*\Cal L_{\Cal E}^{\otimes m}$ is a quasicoherent sheaf on 
$S$ for every $m\geq 0$. Let 
$$
\Gamma_{\Cal E}:=\bigoplus_{m\geq 0} \pi_*\Cal L_{\Cal E}^{\otimes m},
$$
and put 
$$
\tilde{\Sigma}:=\text{Proj}(\Gamma_{\Cal E}).
$$
Since $\Gamma_{\Cal E}$ is a 
quasicoherent sheaf of graded $\Cal O_S$-algebras, 
then $\tilde{\Sigma}$ is a separated scheme over $S$. 

\parindent=12pt

Lest $\Gamma_{\Cal E}$ 
is too big, it can be more convenient to deal with the truncated 
subsheaf of graded $\Cal O_S$-subalgebras,
$$
\Gamma_{\Cal E}(m)\subseteq \Gamma_{\Cal E},
$$
generated by the homogeneous pieces 
$\pi_*\Cal L_{\Cal E}^{\otimes i}$, for $i=0,\dots, m$. We let 
$$
\tilde{\Sigma}(m):=\text{Proj}(\Gamma_{\Cal E}(m))
$$
for every $m>0$. 

Let $\psi\: \Sigma @>>> \tilde{\Sigma}$ 
denote the canonical rational map associated with 
$\Gamma_{\Cal E}$. In addition, we let 
$\psi_m\: \Sigma @>>> \tilde{\Sigma}(m)$ denote the rational map 
associated with $\Gamma_{\Cal E}(m)$, for every $m>0$. 
The scheme $\tilde{\Sigma}$ can be 
regarded as a \sl separated realization of 
$\Sigma$\rm .

\vskip0.4cm

\parindent=0pt

\sl Example 37. \rm We go back to Ex. 34, and use the notations 
therein. We showed that $J_M=J_E^p$ for any non-singular point 
$p\in X$, where $E:=M^*\otimes\Cal O_X(p)$. Even though the scheme 
$J_M$ does not depend on the choice of $p$, the polarization $E$ does. 
Therefore, we can expect the rational map, $\psi\: J_M @>>> \tilde J_M$, 
to depend on $p$. In fact, 
let $L_E$ denote the Theta line bundle on $J_M$ with respect to $E$. 
It follows from the base-change property of the determinant 
of cohomology that $\Gamma_M^*(L_E)\cong\Cal O_X(p)$. So, 
under the identification between $X$ and $J_M$ 
given by $\Gamma_M$, we see that
$\tilde J_M\cong X_p$, where 
$X_p$ is the irreducible component of $X$ containing $p$.

\vskip0.4cm

\parindent=12pt

\proclaim{Theorem 38} Let $\Sigma\subseteq J^{ss}_{\Cal E}$ 
be an open subspace. Let
$$
\psi_m\: \Sigma @>>> \tilde{\Sigma}(m)
$$
denote the canonical rational map, for every $m>0$. If $m>>0$, then 
the following statements hold:
\roster
\item $\psi_m$ is defined on $\Sigma$, and is scheme-theoretically dominant.
\item For every geometric point $s\in S$, each fibre of 
$\psi_m(s)$ is contained in a JH-equivalence class of $\Sigma(s)$.
\item If $\Sigma$ is universally closed over $S$, and $J^s_{\Cal E}\subseteq 
\Sigma$, then the restriction $\left.\psi_m\right|_{J^s_{\Cal E}}$ is 
an open embedding.
\endroster
\endproclaim

\demo{Proof} We show (1) first. 
Let $[I]\in\Sigma$ be a geometric point, representing a semistable 
sheaf on $X(s)$, where $s\in S$ is the geometric point 
below $[I]$. We want to show that $\psi_m$ is defined at $[I]$. 
Since the formation of 
$\Gamma_{\Cal E}$ commutes with \'etale base change, 
we may replace $S$ by any of its \'etale coverings. 
It follows from Thm. 8, Rmk. 15, Rmk. 16 and Rmk. 30 that, 
replacing $S$ by an \'etale covering, if necessary, there is a 
vector bundle $\Cal F$ on $X$, with rank $mr$ and $\det\Cal F\cong
(\det\Cal E)^{\otimes m}$, such that
$$
h^0(X(s),I\otimes\Cal F(s))=h^1(X(s),I\otimes\Cal F(s))=0.\tag{38.1}
$$
We may also assume that $\Cal D_f(\Cal F)\cong\Cal D_f(\Cal E)^{\otimes m}$. 
Let $\theta_{\Cal F}\in H^0(\Sigma,\Cal L^{\otimes m}_{\Cal E})$ denote 
the corresponding theta function. It follows from (38.1) that 
$\theta_{\Cal F}([I])\neq 0$. Thus, $\psi_m$ is defined at $[I]$. It is 
clear that $\psi_m$ is scheme-theoretically dominant.

We prove (2) now. Let $s\in S$ be a geometric point. Let 
$[I],[J]\in\Sigma(s)$ be closed points such that 
$\text{Gr}(I)\not\cong\text{Gr}(J)$. We have to show that 
$\psi_m([I])\neq\psi_m([J])$. We first note that we may replace $S$ by any 
of its \'etale coverings. It follows 
from \cite{\bf 13\rm , Lemma 9 and Rmk. 3} and Rmk. 30 that, 
replacing $S$ by an \'etale covering, if necessary, 
there is a vector bundle $\Cal F$ on $X$, with rank $mr$ and 
$\det\Cal F\cong(\det\Cal E)^{\otimes m}$, such that
$$
h^0(X(s),I\otimes\Cal F(s))=h^1(X(s),I\otimes\Cal F(s))\neq 0\tag{38.2}
$$
and
$$
h^0(X(s),J\otimes\Cal F(s))=h^1(X(s),J\otimes\Cal F(s))=0.\tag{38.3}
$$
We may also assume that $\Cal D_f(\Cal F)\cong\Cal D_f(\Cal E)^{\otimes m}$. 
Let 
$\theta_{\Cal F}\in H^0(\Sigma,\Cal L^{\otimes m}_{\Cal E})$ denote 
the corresponding theta function. It follows from (38.2) and (38.3) 
that $\theta_{\Cal F}([I])=0$, while $\theta_{\Cal F}([J])\neq 0$. 
Thus, $\psi_m([I])\neq\psi_m([J])$.

We prove (3) now. We first note that the statement is local on the \'etale 
topology of $S$. Thus, we may assume that $\Sigma$ is a scheme. 
Since $\Sigma$ is 
universally closed over $S$, and $\psi_m$ is dominant, 
then $\psi_m$ is surjective and closed. It 
follows from the already proved statement (2) that
$$
\psi_m(J^s_{\Cal E})\cap\psi_m(\Sigma\setminus J^s_{\Cal E})=\emptyset.
$$
Hence, $\psi_m(J^s_{\Cal E})\subseteq 
\tilde{\Sigma}(m)$ is open. Moreover, it follows from statement (2) 
that $J^s_{\Cal E}=\psi_m^{-1}(\psi_m(J^s_{\Cal E}))$, and the restriction,
$$
\phi_m:=\left.\psi_m\right|_{J^s_{\Cal E}} \: J^s_{\Cal E} @>>> 
\psi_m(J^s_{\Cal E}),
$$
is bijective. Since $\Sigma$ is universally closed over $S$, then 
$\phi_m$ is universally closed. Since 
$J^s_{\Cal E}$ is separated over $S$, then $\phi_m$ is proper. Since 
$\phi_m$ is bijective and proper, then $\phi_m$ is finite. To finish the 
proof, we need only show that $\phi_m$ separates tangent vectors. Consider 
the natural homomorphism,
$$
h\: \phi_m^*\Omega^1_{\psi_m(J^s_{\Cal E})/S} @>>> \Omega^1_{J^s_{\Cal E}/S}.
$$
Let $\Cal H:=\text{coker}(h)$. We need to show that 
$\Cal H=0$. Suppose, by contradiction, that $\Cal H\neq 0$. Then, 
there are a geometric point $s\in S$, a closed point 
$[I]\in J^s_{\Cal E}(s)$, and a non-zero 
tangent vector $v\in T_{[I],\Sigma(s)}$ such that 
$(d\psi_m)_{[I]}(v)=0$. It follows 
from \cite{\bf 13\rm , Lemma 11 and Rmk. 3} and Rmk. 30 that, 
replacing $S$ by an \'etale covering, if necessary, 
there is a vector bundle $\Cal F$ on $X$, with rank $mr$ and 
$\det\Cal F\cong(\det\Cal E)^{\otimes m}$, such that 
$[I]\in\Theta_{\Cal F}(s)$, but 
$v\not\in T_{[I],\Theta_{\Cal F}(s)}$, where $\Theta_{\Cal F}\subseteq\Sigma$ 
is the zero-scheme of the theta function $\theta_{\Cal F}$ on $\Sigma$. 
Thus, $(d\psi_m)_{[I]}(v)\neq 0$, reaching a contradiction. The 
proof of the theorem is complete.\qed
\enddemo

\vskip0.4cm

\parindent=0pt

\sl Remark 39. \rm Strictly speaking, Lemmas 9 and 11 in \cite{\bf 13} 
hold for curves defined over an algebraically closed field $k$. But, 
applying the same kind of argument used in Rmk. 15, we can show that the 
said lemmas hold for curves defined over any field $k$, up to 
replacing $k$ by one of its separable finite field extensions.

\parindent=12pt

\vskip0.4cm

\proclaim{Corollary 40} $J^s_{\Cal E}$ is a scheme.
\endproclaim

\demo{Proof} Apply statement (3) of Thm. 38 to 
$\Sigma:=J^{ss}_{\Cal E}$.\qed
\enddemo

\parindent=0pt

\vskip0.4cm

\sl Proof of Theorem B. \rm 
By Lemma 29, there are relative polarizations $\Cal E_j$ on 
$X$ over $S$ such that
$$
J=\bigcup_j J^s_{\Cal E_j}.
$$
By Cor. 40, the algebraic space $J^s_{\Cal E_j}$ is a scheme, 
for every $j$. The proof is complete.\qed

\vskip0.4cm

\parindent=0pt

\sl Remark 41. \rm I conjecture that 
$J$ is a scheme if the irreducible components of every fibre of $f$ are 
geometrically integral. 
I can prove my conjecture when $S$ is (the spectrum of) a 
field, by showing that there are enough polarizations on $X$. In this case, 
the argument for constructing polarizations on $X$ is just an adaptation 
of the argument used in the proof of Lemma 29. If one could prove the 
conjecture in general, then one would obtain a true generalization of 
Mumford's result \cite{\bf 9\rm , p. 210}.

\vskip0.4cm

\parindent=12pt

\proclaim{Corollary 42} There is an \'etale covering $S' @>>> S$ such that 
$J\times_S S'$ is a scheme.
\endproclaim

\demo{Proof} Immediate from Lemma 28 and Thm. B.\qed
\enddemo

\vskip0.4cm

\parindent=0pt

\sl Remark 43. \rm If $\pi_*\Cal L_{\Cal E}^{\otimes m}$ is coherent 
for every $m\geq 0$, then $\tilde{\Sigma}(m)$ is locally projective 
over $S$. In this case, if $S$ is excellent, then it follows from 
statement (1) of Thm. 38 that $\tilde{\Sigma}=\tilde{\Sigma}(m)$ 
for $m>>0$. If $\Sigma$ is proper over $S$, then 
$\pi_*\Cal L_{\Cal E}^{\otimes m}$ is coherent for every $m\geq 0$. 
That is the case when $\Sigma=J^{\sigma}_{\Cal E}$. 
We claim that $\pi_*\Cal L_{\Cal E}^{\otimes m}$ is coherent for 
every $m\geq 0$ when $\Sigma=J^{qs}_{\Cal E}$ as well. To show our 
claim, we may replace $S$ by an \'etale covering, if necessary, 
and thus assume that 
there are enough sections $\sigma_1,\dots,\sigma_n\: S @>>> X$ of $f$ through 
the $S$-smooth locus of $X$ such that $J^{qs}_{\Cal E}=J^{\sigma_1}_{\Cal E}
\cup\dots\cup J^{\sigma_n}_{\Cal E}$ (see Lemma 28). Since our claim 
holds for $\Sigma=J^{\sigma_i}_{\Cal E}$, for $i=1,\dots,n$, in the place of 
$J^{qs}_{\Cal E}$, then it holds also for $J^{qs}_{\Cal E}$.

\vskip0.4cm

\parindent=12pt

We observe that, in the proof of Thm. 38, we used 
only theta functions. Hence, we have actually proved a 
more refined result, to be stated more precisely below, in a restricted 
situation. 

\vskip0.4cm

\parindent=0pt

\subheading{5.4. The sheaf of theta functions} 
Assume that $S$ is locally of finite type over an algebraically closed field 
$k$. For every integer $m\geq 0$, and every open subscheme 
$U\subseteq S$, let
$$
A^*_m(U)\subseteq H^0(\Sigma\times_S U,\Cal L_{\Cal E}^{\otimes m})
$$
denote the submodule generated by theta functions, $\theta_{\Cal F}$, 
associated with vector 
bundles $\Cal F$ on $X\times_S U$, of rank $mr$ and 
$\det\Cal F\cong(\det\Cal E)^{\otimes m}\otimes\Cal O_U$, such that 
$\Cal D_{f_U}(\Cal F)\cong\Cal D_f(\Cal E)^{\otimes m}\otimes\Cal O_U$, where 
$f_U\: X\times_S U @>>> U$ is the projection morphism. Of course, 
$A^*_m$ is a presheaf. Let $A_m$ denote the associated 
sheaf in the Zariski topology.

\parindent=12pt

\vskip0.4cm

\proclaim{Proposition 44} For every $m\geq 0$, the sheaf $A_m$ is coherent.
\endproclaim

\demo{Proof} Fix an integer $m\geq 0$. Let 
$\Cal N$ be a relatively ample sheaf on $X$ over $S$. We may assume 
that $S$ is an affine scheme of finite type over $k$. Let $c_0$ be an integer 
such that, for every $c\geq c_0$, and every closed point $s\in S$, the 
sheaf $\Cal E(s)\otimes\Cal N(s)^{\otimes c}$ is generated by 
global sections, and
$$
H^0(X(s),\Cal N^{\otimes -cmr}(s)\otimes(\det\Cal E(s))^{\otimes -m})=0.
\tag{44.1}
$$
We will make use of the following auxiliary sheaves: 
for every integer $c\geq c_0$, and every open subscheme 
$U\subseteq S$, let $A_{m,c}^*(U)\subseteq A_m^*(U)$ denote the submodule 
generated by theta functions associated to vector bundles 
$\Cal F$ on $X\times_S U$, with rank $mr$ and $\Cal D_{f_U}(\Cal F)\cong
\Cal D_f(\Cal E)^{\otimes m}\otimes\Cal O_U$, 
such that $\Cal F$ fits in the middle 
of a short exact sequence of the form:
$$
0 @>>> (\Cal N^{\otimes -c})^{\oplus (mr-1)}\otimes\Cal O_U
@>>> \Cal F @>>> 
(\det\Cal E)^{\otimes m}\otimes\Cal N^{\otimes c(mr-1)}\otimes\Cal O_U
@>>> 0.\tag{44.2}
$$
Let $A_{m,c}$ denote the associated sheaf in the Zariski topology.

Fix $c\geq c_0$. We will first show that $A_{m,c}$ is quasicoherent. 
Replacing $S$ by an affine Zariski covering, if necessary, we may assume 
that $\Cal E$ fits in the middle of an exact sequence of the form:
$$
0 @>>> (\Cal N^{\otimes -c})^{\oplus (mr-1)} 
@>>> \Cal E @>>> 
(\det\Cal E)^{\otimes m}\otimes\Cal N^{\otimes c(mr-1)} @>>> 0.\tag{44.3}
$$
We let $\Cal P_c$ denote the collection of vector bundles 
$\Cal F$ on $X$, with rank $mr$ and 
$\Cal D_f(\Cal F)\cong\Cal D_f(\Cal E)^{\otimes m}$, such that 
$\Cal F$ fits in the middle of a short exact sequence of the form 
(44.2), with $U=S$. Let 
$P_c$ denote the free $\Cal O_S$-module with basis $\Cal P_c$. Let 
$h_c\: P_c @>>> A_{m,c}$ denote the $\Cal O_S$-module homomorphism 
mapping $\Cal F\in\Cal P_c$ to $\theta_{\Cal F}$, 
for every $\Cal F\in\Cal P_c$. 
Since $\pi_*\Cal L_{\Cal E}^{\otimes m}$ is quasicoherent, in order to 
show that $A_{m,c}$ is quasicoherent we need only prove that $h_c$ is 
surjective. Let $s\in S$ be a closed point.  Let 
$U\subseteq S$ be an open neighbourhood of $s$, and let 
$\theta\in A_{m,c}(U)$. Since 
$A_{m,c}$ is the sheafification of $A^*_{m,c}$, replacing $U$ by an 
open neighbourhood of $s$, if necessary, 
we may assume that there are finitely many vector bundles 
$\Cal G_1,\dots,\Cal G_t$ on $X\times_S U$ of rank $mr$, fitting in the 
middle of short exact sequences of the form:
$$
0 @>>> (\Cal N^{\otimes -c})^{\oplus (mr-1)}\otimes\Cal O_U
@>>> \Cal G_i @>>> 
(\det\Cal E)^{\otimes m}\otimes \Cal N^{\otimes c(mr-1)}\otimes\Cal O_U 
@>>> 0,\tag{44.4}
$$
such that
$$
\theta=a_1\theta_{\Cal G_1}+\dots+a_t\theta_{\Cal G_t},
$$
for certain $a_1,\dots,a_t\in H^0(U,\Cal O_U)$. We may also assume 
that $U=S_b$, for some $b\in H^0(S,\Cal O_S)$. Let
$$
W:=H^1(X\times S,(\Cal N^{\otimes -cmr})^{\oplus (mr-1)}\otimes
(\det\Cal E)^{\otimes -m}).
$$
It follows from (44.1) that the formation of $W$ commutes with base change. 
It follows from (44.4) that each $\Cal G_i$ is ``represented'' by an 
element $g_i\in W_b$. Actually, we may assume that each $g_i$ comes 
from an element $f_i\in W$. In other words, for each $i$ there is a 
vector bundle $\Cal F_i$ of rank $mr$ on $X$, fitting in the middle of a short 
exact sequence of the form:
$$
0 @>>> (\Cal N^{\otimes -c})^{\oplus (mr-1)} @>>> \Cal F_i @>>> 
(\det\Cal E)^{\otimes m}\otimes \Cal N^{\otimes c(mr-1)} @>>> 0,\tag{44.5}
$$
such that the restriction of (44.5) to $X\times_S U$ is (44.4). It follows 
from (44.3) and (44.5), and the additive property of the determinant of 
cohomology, that $\Cal D_f(\Cal F_i)\cong\Cal D_f(\Cal E)^{\otimes m}$. 
Thus, $\Cal F_i\in\Cal P_c$ for $i=1,\dots,t$. Since 
$\Cal G_i=\left.\Cal F_i\right|_{X\times_S U}$, then 
$\theta_{\Cal G_i}=\left.\theta_{\Cal F_i}\right|_{\Sigma\times_S U}$ 
(modulo $H^0(U,\Cal O_U)^*$) for $i=1,\dots,t$. So, $h_c$ is surjective, 
completing the proof that $A_{m,c}$ is quasicoherent.

Since $A_m$ is generated, as a subsheaf of $\pi_*\Cal L_{\Cal E}^{\otimes m}$, 
by the subsheaves $A_{m,c}$, for $c\geq c_0$, then 
it follows that $A_m$ is quasicoherent. 
We will now show that $A_m$ is coherent. 
Since $A_m$ is quasicoherent, and the formation of 
$\pi_*\Cal L_{\Cal E}^{\otimes m}$ commutes with \'etale base change, 
then we may replace $S$ by any of its \'etale coverings. We may thus 
assume (see Lemmas 28 and 29) that there are a 
section $\sigma\: S @>>> X$ of $f$ through the $S$-smooth 
locus of $X$, and finitely many relative polarizations 
$\Cal H_1,\dots,\Cal H_t$ on $X$ over $S$, such that
$$
\Sigma\subseteq \Sigma':=\bigcup_{j=1}^t J^{\sigma}_{\Cal H_i}.
$$
By Thm. A, the spaces $J^{\sigma}_{\Cal H_i}$ are 
proper over $S$. If $\Sigma=\Sigma'$, then $\pi_*\Cal L_{\Cal E}^{\otimes m}$ 
is coherent, and thus so is $A_m$. At any rate, since theta functions 
over $\Sigma$ are restrictions of 
theta functions over $\Sigma'$, it follows that $A_m$ is 
coherent. The proof is complete.\qed
\enddemo

\vskip0.4cm

Let
$$
V_{\Cal E}:=\bigoplus_{m\geq 0} A_m\subseteq\Gamma_{\Cal E}.
$$
It follows from the additive property of the determinant of cohomology, and 
the additive property of $\Phi_{\Cal F/\Cal E}$ 
\cite{\bf 14\rm , Thm. I.1, p. 509}, that 
$V_{\Cal E}$ is a graded subsheaf of $\Cal O_S$-subalgebras of 
$\Gamma_{\Cal E}$. The sheaf $V_{\Cal E}$ is called the \sl sheaf of theta 
functions \rm on $\Sigma$ with respect to $\Cal E$. Put
$$
\overline{\Sigma}:=\text{Proj}(V_{\Cal E}).
$$
It follows from Prop. 44 that 
$\overline{\Sigma}$ is separated over $S$. For every 
$m>0$, denote by $V_{\Cal E}(m)$ the graded subsheaf of 
$\Cal O_S$-subalgebras of $V_{\Cal E}$ 
generated by the homogeneous pieces of $V_{\Cal E}$ of degree at most 
$m$. Let
$$
\overline{\Sigma}(m):=\text{Proj}(V_{\Cal E}(m))
$$
for every $m>0$. It follows from Prop. 44 that 
$\overline{\Sigma}(m)$ is locally projective over $S$ for every $m>0$.

\vskip0.4cm

\parindent=0pt

\sl Proof of Theorem C. \rm It is a simple observation that the 
same statements we proved in Thm. 38 remain valid if we replace 
$\tilde{\Sigma}(m)$ by $\overline{\Sigma}(m)$. Since statement (1) 
of Thm. 38 holds for $\overline{\Sigma}(m)$ for $m>>0$, and 
$A_m$ is coherent for every $m$, then it follows that 
$\overline{\Sigma}=\overline{\Sigma}(m)$ for $m>>0$. In particular, 
$\overline{\Sigma}$ is locally projective over $S$. 
The proof is complete.\qed

\vskip0.4cm

\parindent=12pt

Note that, if $\Sigma'\subseteq\Sigma$ is an open subspace, then 
there is a natural closed embedding, 
$\overline{\Sigma}'\subseteq\overline{\Sigma}$.

\vskip0.8cm

\parindent=12pt

\heading
6. Comparison with Seshadri's compactification
\endheading

Let $X$ be a reduced curve over an algebraically closed field $k$. 
Let $X_1,\dots,X_n$ denote the irreducible components of $X$. 
Let $\underline a:=(a_1,\dots,a_n)$ 
be a $n$-uple of positive rational numbers such that 
$a_1+\dots+a_n=1$. According to Seshadri's definition \cite{\bf 27\rm , 
Part 7}, a 
torsion-free, rank 1 sheaf $I$ on $X$ is $\underline a$-semistable 
(resp. $\underline a$-stable) if and only if
$$
\chi(I_Y)\geq (a_{i_1}+\dots+a_{i_m})\chi(I) 
\text{\  \  (resp. \  } 
\chi(I_Y)>(a_{i_1}+\dots+a_{i_m})\chi(I)
\text{)}
$$
for every non-empty, proper subcurve $Y=X_{i_1}\cup\dots\cup X_{i_m}
\subsetneqq X$.

\parindent=0pt

\vskip0.4cm

\sl Remark 45. \rm 
We shall observe that Seshadri's notions of semistability and stability 
are equivalent to ours. 
Let $d$ be a positive integer. 
Let $\underline a$ be a $n$-uple of 
positive rational numbers such that $a_1+\dots +a_n=1$. 
Let $\underline A:=(A_1,\dots, A_n)$ be a $n$-uple of 
positive integers such that
$$
A_ia_j=A_ja_i \text{\  for \  } 1\leq i,j \leq n;\tag{45.1}
$$
and $A:=A_1+\dots+A_n$ is a multiple of $d$. 
Since $a_1+\dots+a_n=1$, then it follows from (45.1) that $a_i:=A_i/A$ 
for $i=1,\dots,n$. For each $i=1,\dots,n$, 
let $x^i_1,\dots,x^i_{A_i}\in X_i$ be 
distinct, non-singular points of $X$. Let
$$
\Cal O_X(1):=\Cal O_X(\sum \Sb 1\leq i\leq n\\ 1\leq j\leq A_i \endSb 
x^i_j).
$$
Of course, $\Cal O_X(1)$ is an ample sheaf on $X$. Let
$$
E:=(\Cal O_X^{\oplus (r/t-1)}\oplus\Cal O_X(-1))^{\oplus t}\tag{45.2}
$$
for a certain integer $t>0$, to be specified later, and $r:=tA/d$. 
Then, the reader can easily verify (see \cite{\bf 13\rm , Rmk. 12}) 
that a torsion-free, rank 1 
sheaf $I$ on $X$ with $\chi(I)=d$ is $\underline a$-semistable 
(resp. $\underline a$-stable) if and only if $I$ is semistable 
(resp. stable) with respect to $E$.

\vskip0.4cm

\parindent=12pt

Let $S(\underline a,d)$ 
(resp. $S'(\underline a,d)$) 
denote the set of isomorphism
classes of $\underline a$-semistable (resp. $\underline a$-stable) 
torsion-free, rank 1 sheaves of Euler characteristic $d$ on $X$. 

\vskip0.4cm

\proclaim{Theorem 46} 
{\rm (Seshadri)} 
There is a coarse moduli space
for $S'(\underline a,d)$, whose underlying scheme is a
quasi-projective variety denoted by 
$U^s(\underline a,d)$. Moreover, 
$U^s(\underline a,d)$ 
has a natural projective compactification, 
denoted by $U(\underline a,d)$. 
The set of closed points of $U(\underline a,d)$ is
isomorphic to the quotient of $S(\underline a,d)$ 
by the JH-equivalence relation.\qed
\endproclaim

\demo{Proof} See \cite{\bf 27\rm , Thm. 15, p. 155}.\qed
\enddemo

\vskip0.4cm

$U(\underline a,d)$ was obtained by Seshadri as a good 
quotient, $\phi\: R @>>> U(\underline a,d)$, 
of a quasi-projective scheme $R$ under the action of a 
reductive group, using Geometric Invariant Theory. 
In Seshadri's construction, there is a relatively 
torsion-free, rank 1 sheaf $\text{\bf I}$ on $X\times R$ such that 
$\text{\bf I}$ has the local universal property for $S(\underline a,d)$. 
If we take the polarization $E$ in (45.2) with $t>>0$, 
then the determinant of cohomology, 
$\Cal D(\text{\bf I}\otimes E)$ (with respect to the projection 
$X\times R @>>> R$), descends to an 
ample invertible sheaf $\overline L_E$ on $U(\underline a,d)$ 
(see \cite{\bf 13\rm , Thm. 13} and the argument thereafter). 
In addition, the global sections $\sigma_{\text{\bf I}\otimes F}$ of 
$\Cal D(\text{\bf I}\otimes F)$, corresponding to vector bundles $F$ on $X$, 
with rank $mr$ and 
$\det F\cong(\det E)^{\otimes m}$ for some $m>0$, induce in a natural way 
(modulo $k^*$) 
global sections of $\Cal D(\text{\bf I}\otimes E)^{\otimes m}$ that 
descend to sections of $\overline L_E^{\otimes m}$. We shall 
denote by $\overline{\theta}_F$ the global section of 
$\overline L_E^{\otimes m}$ naturally associated with 
$\sigma_{\text{\bf I}\otimes F}$. We 
emphasize that $\overline{\theta}_F$ is well defined modulo $k^*$. Let
$$
\overline{\Gamma}_E:=\bigoplus_{m\geq 0} H^0(U(\underline a,d),
\overline L_E^{\otimes m}),
$$
and let $\overline V_E\subseteq \overline{\Gamma}_E$ 
denote the graded $k$-subalgebra generated by the sections 
$\overline{\theta}_F$. Let
$$
\overline U(\underline a,d):=\text{Proj}(\overline V_E).
$$

\vskip0.4cm

\proclaim{Theorem 47} The rational map,
$$
\pi\: U(\underline a,d) @>>> \overline U(\underline a,d),
$$
is defined everywhere, bijective, and an isomorphism on $U^s(\underline a,d)$.
\endproclaim

\demo{Proof} See \cite{\bf 13\rm , Thm. 15 and Thm. 16}.\qed
\enddemo

\vskip0.4cm

Let $J^{ss}_E$ denote the fine moduli space of semistable sheaves on $X$ with 
respect to $E$ (see Sect. 3). Let
$$
\psi\: J^{ss}_E @>>> \overline J^{ss}_E
$$
denote the naturally defined rational map (see Sect. 5). 
Since $J^{ss}_E$ is a 
fine moduli space, there is a naturally defined morphism,
$$
\rho\: J^{ss}_E @>>> U(\underline a,d),
$$
mapping $[I]\in J^{ss}_E$ to $[\text{Gr}(I)]\in U(\underline a,d)$. 
It is clear that $\rho^{-1}(U^s(\underline a,d))=J^s_E$, and 
that $\left.\rho\right|_{J^s_E}\: J^s_E @>>> U^s(\underline a,d)$ 
is an isomorphism. 
Let $L_E$ denote the Theta line bundle on 
$J^{ss}_E$. Since the relatively torsion-free, rank 1 sheaf $\text{\bf I}$ on 
$X\times R$ has the local universal property for $S(\underline a,d)$, 
then the morphism $\rho$ factors locally through the quotient 
morphism $\phi\: R @>>> U(\underline a,d)$. Thus, 
$L_E\cong\rho^*\overline L_E$ locally on $J^{ss}_E$. Actually, 
it is easy to show that the local isomorphisms patch to a global isomorphism,
$$
L_E\cong\rho^*\overline L_E.\tag{47.1}
$$
We can also show that $\rho^*\overline{\theta}_F=\theta_F$ 
(modulo $k^*$) under the above isomorphism. Let 
$\Sigma\subseteq J^{ss}_{\Cal E}$ be an open subscheme. 
Let $V_E$ denote the ring of theta functions on $\Sigma$ (see Subsect. 5.4). 
The isomorphism 
(47.1) induces a surjective graded homomorphism of $k$-algebras,
$$
\overline V_E @>>> V_E.
$$
The above surjective homomorphism induces a closed embedding,
$$
\iota\: \overline{\Sigma} \hookrightarrow \overline U(\underline a,d).
$$
It is clear that
$$
\pi\circ\left.\rho\right|_{\Sigma}=\iota\circ\psi.
$$

If $\Sigma\supseteq J^p_E$ for a certain non-singular $p\in X$, 
then it follows from Thm. 7 that $\left.\rho\right|_{\Sigma}$ 
is surjective. Since $\pi$ is bijective, it follows that $\iota$ 
is also bijective. 
If $U(\underline a,d)$ were reduced, then $\iota$ would be an 
isomorphism. If $X$ has at most ordinary double points for singularities, then 
$U(\underline a,d)$ is known to be reduced \cite{\bf 1\rm , Cor. 3.5}.

\vskip0.4cm

\proclaim{Theorem 48} If $X$ has at most ordinary double points for 
singularities, then
$$
\iota\: \overline{\Sigma}\cong \overline U(\underline a,d)
$$
is an isomorphism.\qed
\endproclaim

\vskip0.4cm

If one had a better handling on the tangent spaces at semistable points 
of $U(\underline a,d)$, then one could say more about the morphisms 
$\pi$ and $\iota$. 
Unfortunately though, Geometric Invariant Theory does not 
seem to provide good infinitesimal information about the quotients obtained 
by its method, in general. 

\vskip0.8cm

\heading
Acknowledgements
\endheading

I take the opportunity to thank my Ph.D. advisor, 
Prof. Kleiman, for introducing me to the 
compactification problem and for several comments on different versions 
of the present article. I am also grateful to Prof. Ramero for 
telling me about Faltings' construction of the moduli space of  
semistable vector bundles on a smooth curve via theta functions. I 
would like to thank Waseda University, specially Prof. Kaji, 
Prof. Morimoto, Prof. Ohno and Prof. Hara, for the warm hospitality extended 
during the period the essential part of this work was finished. 
I am also grateful to Prof. Homma for all the help received from him 
during the said period.

\vskip0.4cm

\parindent=12pt

\eightsmc Instituto de Matem\'atica Pura e Aplicada, 
Estrada Dona Castorina 110, 22460-320 Rio de Janeiro RJ, Brazil.

\smallslant E-mail address: \smallrm esteves\@ impa.br


\vskip0.8cm

\Refs

\ref \no 1 \by V. Alexeev \paper Compactified Jacobians 
\jour E-print at alg-geom\@ eprints.math.duke.edu/9608012, August 1996
\endref

\ref \no 2 \by A. Altman, A. Iarrobino and S. Kleiman \paper 
Irreducibility of the compactified Jacobian \jour Real and complex 
singularities, Oslo 1976 (Proc. 9th nordic summer school), 
Sijthoff and Noordhoff, 1977 \pages pp. 1--12
\endref

\ref \no 3 \by A. Altman and S. Kleiman \paper 
Compactifying the Picard scheme II 
\jour Amer. J. Math. \vol 101 \yr 1979 \pages 10--41
\endref

\ref \no 4 \by A. Altman and S. Kleiman \paper 
Compactifying the Picard scheme \jour Adv. Math. \vol 35 \yr 1980 
\pages 50--112
\endref

\ref \no 5 \by A. Altman and S. Kleiman 
\paper Algebraic systems of 
divisor-like subschemes 
\jour Compositio Math. \vol 29 \yr 1974 \pages 113--139
\endref

\ref \no 6 \by E. Arbarello, M. Cornalba, P. Griffiths and J. 
Harris 
\paper Geometry of algebraic curves, Vol. I \jour Grundlehren der 
mathematischen Wissenschaften, vol. 267, Springer-Verlag, 1985
\endref

\ref \no 7 \by M. Artin \paper Algebraization of formal moduli I 
\jour Global Analysis, Papers in honor of K. Kodaira, University 
of Tokyo Press, Princeton University Press, 1969 \pages pp. 21--71
\endref

\ref \no 8 \by A. Beauville \paper Counting rational curves on 
K3 surfaces \jour Available from the e-print service 
at alg-geom\@ eprints.math.duke.edu/9701019, 
January 1997
\endref

\ref \no 9 \by S. Bosch, W. L\"utkebohmert, and M. Raynaud \paper 
N\'eron models \jour Ergebnisse der Mathematik 
und ihrer Grenzgebiete, vol. 21, Springer Verlag, 1990
\endref

\ref \no 10 \by L. Caporaso \paper A compactification of the universal 
Picard variety over the moduli space of stable curves 
\jour J. Amer. Math. Soc. \vol 7 \yr 1994 \pages 589--660
\endref

\ref \no 11 \by C. D'Souza \paper Compactification of generalized 
Jacobians \jour Proc. Indian Acad. Sci. Sect. A Math. Sci. 
\vol 88 \yr 1979 \pages 419--457
\endref

\ref \no 12 \by E. Esteves \paper Very ampleness for Theta on the 
compactified Jacobian \jour Available from the e-print service 
at alg-geom\@ eprints.math.duke.edu/9709005, September 1997. 
To appear in Math. Z.
\endref

\ref \no 13 \by E. Esteves \paper Separation properties of 
theta functions \jour Available from the e-print service 
at alg-geom\@ eprints.math.duke.edu, 
September 1997
\endref

\ref \no 14 \by G. Faltings \paper Stable {\eightsl G}-bundles and projective 
connections \jour J. Algebraic Geom. \vol 2 \yr 1993 
\pages 507--568
\endref

\ref \no 15 \by D. Gieseker \paper Moduli of curves \jour Tata Inst. Fund. 
Res. Lecture Notes, Springer-Verlag, 1982
\endref

\ref \no 16 \by A. Grothendieck with J. Dieudonn\'e \paper 
\'Elem\'ents de G\'eom\'etrie Alg\'ebrique 
III-1, IV-4 \jour Inst. Hautes \'Etudes Sci. Publ. Math. 
\vol 11\eightrm , 
\eightbf 32 \yr 1961, 1967
\endref

\ref \no 17 \by A. Grothendieck \paper Fondements de la g\'eometrie 
alg\'ebrique \jour Semin\'aire Bourbaki, exp. \eightbf 232 
\yr 1961/62
\endref

\ref \no 18 \by J. Igusa \paper Fiber systems of Jacobian varieties 
\jour Amer. J. Math. \vol 78 \yr 1956 \pages 171--199
\endref

\ref \no 19 \by M. Ishida \paper Compactifications of a family of 
generalized Jacobian varieties \jour Proc. Internat. Sympos. Algebraic 
Geometry, Kyoto, January, 1977 (M. Nagata, ed.), 
Kinokuniya, Tokyo, 1978 \pages pp. 503--524
\endref

\ref \no 20 \by F. Knudsen and D. Mumford 
\paper The projectivity of the moduli 
space of stable curves I: Preliminaries on ``det'' and ``Div'' 
\jour Math. Scand. \vol 39 \yr 1976 
\pages 19--55
\endref

\ref \no 21 \by S. Langton \paper Valuative criteria for 
families of vector bundles on algebraic varieties \jour 
Ann. Math. \vol 101 \yr 1975 \pages 88--110
\endref

\ref \no 22 \by A. Meyer and D. Mumford \paper Further comments on 
boundary points \jour Unpublished lecture notes distributed at the 
Amer. Math. Soc. Summer Institute, Woods Hole, 1964
\endref

\ref \no 23 \by D. Mumford \paper Lectures on curves on an 
algebraic surface \jour Annals of Math. Studies, no. 59, 
Princeton University Press, 1966
\endref

\ref \no 24 \by T. Oda and C.S. Seshadri \paper Compactifications of the 
generalized Jacobian variety \jour Trans. Amer. Math. Soc. \vol 253 
\yr 1979 \pages 1--90
\endref

\ref \no 25 \by R. Pandharipande \paper 
A compactification over $\overline{\text{\eightsl M}}_{\text{\smallit g}}$ 
of the universal moduli space of slope-semistable vector bundles 
\jour J. Amer. Math. Soc. \vol 9 \yr 1996 \pages 425--471
\endref

\ref \no 26 \by M. Raynaud (written by S. Kleiman) \paper 
Un th\'eor\`eme de representabilit\'e relative sur le foncteur de Picard 
\jour Lecture Notes in Mathematics, vol. 225, Springer-Verlag, 1971 
\pages pp. 595--615
\endref

\ref \no 27 \by C.S. Seshadri \paper Fibr\'es vectoriels
sur les courbes alg\'ebriques \jour Ast\'erisque \vol 96 \yr 1982
\endref

\ref \no 28 \by C.S. Seshadri \paper Vector bundles on
curves \jour Contemp. Math. \vol 153 \yr 1993 \pages 163--200
\endref

\ref \no 29 \by C. Simpson \paper Moduli of representations of the 
fundamental group of a smooth projective variety I 
\jour Inst. Hautes \'Etudes Sci. Publ. Math. 
\vol 79 \yr 1994 \pages 47--129
\endref

\endRefs
\enddocument